\documentclass[10pt]{amsart}

\newcommand{\bff}{\bf}

\usepackage{graphics}
\usepackage{amssymb}
\usepackage{amsxtra}

\input diagrams

\usepackage[arrow, matrix]{xy}

\usepackage{ifthen}

\theoremstyle{plain}

\newtheorem{theo}{Theorem}[section]
\newtheorem{lemm}[theo]{Lemma}
\newtheorem{prop}[theo]{Proposition}
\newtheorem{coro}[theo]{Corollary}

\theoremstyle{definition}

\newtheorem{defi}[theo]{Definition}
\newtheorem{rema}[theo]{Remark}

\newfont{\rmm}{cmr10 scaled 1000}
\newfont{\itt}{cmsl10 scaled 1000}

\newfont{\rM}{cmr10 scaled 1700}

\newcommand{\lb}{\label}
\newcommand{\mlb}{\label}

\newcommand{\rrf}[1]{(\ref{#1})}
\newcommand{\mrf}[1]{\ref{#1}}

\newarrow{dashto}{}{dash}{}{dash}>

\newarrow{mapsto}|--->

\newarrow{dotsto}|...>

\begin{document}

%ii\input /home/a/auxxx.doc
%%%%%%%%%%%%%%%%%
%%%%ALFAVITY
%%%%%%%%%%%%%%%%%

\renewcommand{\a}{\alpha}
\renewcommand{\b}{\beta}
\newcommand{\g}{\gamma}
\renewcommand{\d}{\delta}
\newcommand{\e}{\epsilon}
\newcommand{\ve}{\varepsilon}
\newcommand{\z}{\zeta}
\renewcommand{\t}{\theta}
\renewcommand{\l}{\lambda}
\renewcommand{\k}{\varkappa}
\newcommand{\m}{\mu}
\newcommand{\n}{\nu}
\renewcommand{\r}{\rho}
\newcommand{\vr}{\varrho}
\newcommand{\s}{\sigma}
\newcommand{\vp}{\varphi}
\renewcommand{\o}{\omega}
\renewcommand{\Re}{\text{\rm Re }}
%GrekBig

\newcommand{\G}{\Gamma}
\newcommand{\D}{\Delta}
\newcommand{\T}{\Theta}
\renewcommand{\L}{\Lambda}
\renewcommand{\P}{\Pi}
\newcommand{\Si}{\Sigma}
\renewcommand{\O}{\Omega}
\newcommand{\Up}{\Upsilon}

\renewcommand{\AA}{{\mathcal A}}
\newcommand{\BB}{{\mathcal B}}
\newcommand{\CC}{{\mathcal C}}
\newcommand{\DD}{{\mathcal D}}
\newcommand{\EE}{{\mathcal E}}
\newcommand{\FF}{{\mathcal F}}
\newcommand{\GG}{{\mathcal G}}
\newcommand{\HH}{{\mathcal H}}
\newcommand{\II}{{\mathcal I}}
\newcommand{\JJ}{{\mathcal J}}
\newcommand{\KK}{{\mathcal K}}
\newcommand{\LL}{{\mathcal L}}
\newcommand{\MM}{{\mathcal M}}
\newcommand{\NN}{{\mathcal N}}
\newcommand{\OO}{{\mathcal O}}
\newcommand{\PP}{{\mathcal P}}
\newcommand{\QQ}{{\mathcal Q}}
\newcommand{\RR}{{\mathcal R}}
\renewcommand{\SS}{{\mathcal S}}
\newcommand{\TT}{{\mathcal T}}
\newcommand{\UU}{{\mathcal U}}
\newcommand{\VV}{{\mathcal V}}
\newcommand{\WW}{{\mathcal W}}
\newcommand{\XX}{{\mathcal X}}
\newcommand{\YY}{{\mathcal Y}}
\newcommand{\ZZ}{{\mathcal Z}}

\renewcommand{\aa}{{\mathbb{A}}}
\newcommand{\bb}{{\mathbb{B}}}
\newcommand{\cc}{{\mathbb{C}}}
\newcommand{\dd}{{\mathbb{D}}}
\newcommand{\ee}{{\mathbb{E}}}
\newcommand{\ff}{{\mathbb{F}}}
\renewcommand{\gg}{{\mathbb{G}}}
\newcommand{\hh}{{\mathbb{H}}}
\newcommand{\ii}{{\mathbb{I}}}
\newcommand{\jj}{{\mathbb{J}}}
\newcommand{\kk}{{\mathbb{K}}}
\renewcommand{\ll}{{\mathbb{L}}}
\newcommand{\mm}{{\mathbb{M}}}
\newcommand{\nn}{{\mathbb{N}}}
\newcommand{\oo}{{\mathbb{O}}}
\newcommand{\pp}{{\mathbb{P}}}
\newcommand{\qq}{{\mathbb{Q}}}
\newcommand{\rr}{{\mathbb{R}}}
\renewcommand{\ss}{{\mathbb{S}}}
\newcommand{\ttt}{{\mathbb{T}}}
\newcommand{\uu}{{\mathbb{U}}}
\newcommand{\vv}{{\mathbb{V}}}
\newcommand{\ww}{{\mathbb{W}}}
\newcommand{\xx}{{\mathbb{X}}}
\newcommand{\yy}{{\mathbb{Y}}}
\newcommand{\zz}{{\mathbb{Z}}}

\newcommand{\AAA}{{\mathbf{A}}}
\newcommand{\BBB}{{\mathbf{B}} }
\newcommand{\CCC}{{\mathbf{C}} }
\newcommand{\DDD}{{\mathbf{D}} }
\newcommand{\EEE}{{\mathbf{E}} }
\newcommand{\FFF}{{\mathbf{F}} }
\newcommand{\GGG}{{\mathbf{G}}}
\newcommand{\HHH}{{\mathbf{H}}}
\newcommand{\III}{{\mathbf{I}}}
\newcommand{\JJJ}{{\mathbf{J}}}
\newcommand{\KKK}{{\mathbf{K}}}
\newcommand{\LLL}{{\mathbf{L}}}
\newcommand{\MMM}{{\mathbf{M}}}
\newcommand{\NNN}{{\mathbf{N}}}
\newcommand{\OOO}{{\mathbf{O}}}
\newcommand{\PPP}{{\mathbf{P}}}
\newcommand{\QQQ}{{\mathbf{Q}}}
\newcommand{\RRR}{{\mathbf{R}}}
\newcommand{\SSS}{{\mathbf{S}}}
\newcommand{\TTT}{{\mathbf{T}}}
\newcommand{\UUU}{{\mathbf{U}}}
\newcommand{\VVV}{{\mathbf{V}}}
\newcommand{\WWW}{{\mathbf{W}}}
\newcommand{\XXX}{{\mathbf{X}}}
\newcommand{\YYY}{{\mathbf{Y}}}
\newcommand{\ZZZ}{{\mathbf{Z}}}

\newcommand{\gA}{{\mathfrak{A}}}
\newcommand{\gB}{{\mathfrak{B}}}
\newcommand{\gC}{{\mathfrak{C}}}
\newcommand{\gD}{{\mathfrak{D}}}
\newcommand{\gE}{{\mathfrak{E}}}
\newcommand{\gF}{{\mathfrak{F}}}
\newcommand{\gG}{{\mathfrak{G}}}
\newcommand{\gH}{{\mathfrak{H}}}
\newcommand{\gI}{{\mathfrak{I}}}
\newcommand{\gJ}{{\mathfrak{J}}}
\newcommand{\gK}{{\mathfrak{K}}}
\newcommand{\gL}{{\mathfrak{L}}}
\newcommand{\gM}{{\mathfrak{M}}}
\newcommand{\gN}{{\mathfrak{N}}}
\newcommand{\gO}{{\mathfrak{O}}}
\newcommand{\gP}{{\mathfrak{P}}}
\newcommand{\gQ}{{\mathfrak{Q}}}
\newcommand{\gR}{{\mathfrak{R}}}
\newcommand{\gS}{{\mathfrak{S}}}
\newcommand{\gT}{{\mathfrak{T}}}
\newcommand{\gU}{{\mathfrak{U}}}
\newcommand{\gV}{{\mathfrak{V}}}
\newcommand{\gW}{{\mathfrak{W}}}
\newcommand{\gX}{{\mathfrak{X}}}
\newcommand{\gY}{{\mathfrak{Y}}}
\newcommand{\gZ}{{\mathfrak{Z}}}

\newcommand{\gota}{{\mathfrak{a}}}
\newcommand{\gotb}{{\mathfrak{b}}}
\newcommand{\gotc}{{\mathfrak{c}}}
\newcommand{\gotd}{{\mathfrak{d}}}
\newcommand{\gote}{{\mathfrak{e}}}
\newcommand{\gotf}{{\mathfrak{f}}}
\newcommand{\gotg}{{\mathfrak{g}}}
\newcommand{\goth}{{\mathfrak{h}}}
\newcommand{\goti}{{\mathfrak{i}}}
\newcommand{\gotj}{{\mathfrak{j}}}
\newcommand{\gotk}{{\mathfrak{k}}}
\newcommand{\gotl}{{\mathfrak{l}}}
\newcommand{\gotm}{{\mathfrak{m}}}
\newcommand{\gotn}{{\mathfrak{n}}}
\newcommand{\goto}{{\mathfrak{o}}}
\newcommand{\gotp}{{\mathfrak{p}}}
\newcommand{\gotq}{{\mathfrak{q}}}
\newcommand{\gotr}{{\mathfrak{r}}}
\newcommand{\gots}{{\mathfrak{s}}}
\newcommand{\gott}{{\mathfrak{t}}}
\newcommand{\gotu}{{\mathfrak{u}}}
\newcommand{\gotv}{{\mathfrak{v}}}
\newcommand{\gotw}{{\mathfrak{w}}}
\newcommand{\gotx}{{\mathfrak{x}}}
\newcommand{\goty}{{\mathfrak{y}}}
\newcommand{\gotz}{{\mathfrak{z}}}

%%%%%%%%%%%%%%%%%%%%%%%%%%%%%%%%%%%%%%%%%%%%%%%%%
%-----------------------------------------------
%%%%%%%%%%%%%%%%%%%%%%%%%%%%%%%%%%%%%%%%%%%%%%%%%
%-----------------------------------------------
%%%%%%%%%%%%%%%%%%%%%%%%%%%%%%%%%%%%%%%%%%%%%%%%%
%-----------------------------------------------
%%%%%%%%%%%%%%%%%%%%%%%%%%%%%%%%%%%%%%%%%%%%%%%%%
%-----------------------------------------------

%MACROS

\newcommand{\kkrest}{\begin{picture}(14,14)
\put(00,04){\line(1,0){14}}
\put(00,02){\line(1,0){14}}
\put(06,-4){\line(0,1){14}}
\put(08,-4){\line(0,1){14}}
\end{picture}     }

\newcommand{\krest}{~\kkrest~}

%%%%%%%%%%%%%%%%%%%%%%%%%%%
%TEXT/ROMAN
%%%%%%%%%%%%%%%%%%%%%%%%%%%%%

\newcommand{\grd}{{\text{\rm grd}}}
\newcommand{\id}{\text{id}}
\newcommand{\Tb}{\text{ \rm Tb}}
\newcommand{\Log}{\text{\rm Log }}
\newcommand{\Wh}{\text{\rm Wh }}
\newcommand{\Ker}{\text{\rm Ker }}
\newcommand{\Ext}{\text{\rm Ext}}
\newcommand{\Hom}{\text{\rm Hom}}
\newcommand{\diam}{\text{\rm diam}}
\newcommand{\Homb}{\text{\rm Hom}b}
\newcommand{\Lg}{\text{\rm Lg }}
\newcommand{\ind}{\text{\rm ind}}
\newcommand{\rk}{\text{\rm rk }}
\renewcommand{\Im}{\text{\rm Im }}
\newcommand{\supp}{\text{\rm supp }}
\newcommand{\Int}{\text{\rm Int }}
\newcommand{\grad}{\text{\rm grad}}
\newcommand{\Fix}{\text{\rm Fix}}
\newcommand{\Exp}{\text{\rm Exp}}
\newcommand{\Per}{\text{\rm Per}}
\newcommand{\TL}{\text{\rm TL}}
\newcommand{\Id}{\text{\rm Id}}
\newcommand{\Vect}{\text{\rm Vect}}
\newcommand{\vvol}{\text{\rm vol}}
\newcommand{\Mat}{\text\rm Mat}
\newcommand{\Tub}{\text{\rm Tub}}
\newcommand{\Imm}{\text{\rm Im}}
\newcommand{\tn}{\text{\rm t.n.}}
\newcommand{\card}{\text{\rm card }}
\newcommand{\GL}{\text{\rm GL }}

\newcommand{\track}{\text{\rm Track}}
\newcommand{\sgn}{\text{\rm sgn}}
\newcommand{\Arctg}{\text{\rm Arctg }}
\newcommand{\tg}{\text{\rm tg }}
\newcommand{\Arcsin}{\text{\rm Arcsin }}

%%%%%%%%%%%%%%%%%%%%%%%%%%%%%%%%%%
%%%%%%%%%%%%%%%%%%%%%%%%%%%%%%%%%%%%%%%%%%%%
%BEGIN/END
%%%%%%%%%%%%%%%%%%%%%%%%%%%%%%%%%%%%%%%%%%%%%
%%%%%%%%%%%%%%%%%%%%%%%%%%%%%%%%%%%%

\newcommand{\bere}{\begin{rema}}
\newcommand{\bede}{\begin{defi}}

\renewcommand{\beth}{\begin{theo}}
\newcommand{\bele}{\begin{lemm}}
\newcommand{\bepr}{\begin{prop}}
\newcommand{\beeq}{\begin{equation}}
\newcommand{\bega}{\begin{gather}}
\newcommand{\been}{\begin{enumerate}}

\newcommand{\bedee}{\begin{defii}}
\newcommand{\bethh}{\begin{theoo}}
\newcommand{\belee}{\begin{lemmm}}
\newcommand{\beprr}{\begin{propp}}

\newcommand{\beco}{\begin{coro}}

\newcommand{\beal}{\begin{aligned}}

\newcommand{\enre}{\end{rema}}

\newcommand{\enco}{\end{coro}}
\newcommand{\enpr}{\end{prop}}
\newcommand{\enth}{\end{theo}}
\newcommand{\enle}{\end{lemm}}
\newcommand{\enen}{\end{enumerate}}
\newcommand{\enga}{\end{gather}}
\newcommand{\eneq}{\end{equation}}
\newcommand{\enal}{\end{aligned}}

\newcommand{\bq}{\begin{equation}}
\newcommand{\bqq}{\begin{equation*}}

%%%%%%%%%%%%%%%%%%%
%%%%%%%%%%%%%%%%%%%%
%%%DIVERSE ABBREVIATIONS FOR TEX MACROS
%%%%%%%%%%%%%%%%%%%
%%%%%%%%%%%%%%%%%%%%

\renewcommand{\leq}{\leqslant}
\renewcommand{\geq}{\geqslant}
\newcommand{\vphi}{\varphi}
\newcommand{\vide}{  \emptyset  }
\newcommand{\bu}{\bullet}
\newcommand{\pfff}{\pitchfork}
\newcommand{\mx}{\mbox}

\newcommand{\mxx}[1]{\quad\mbox{#1}\quad}
 \newcommand{\mxxx}[1]{\hspace{0.1cm}\mbox{#1} \quad  }
\newcommand{\wi}{\widetilde}

\newcommand{\ove}{\overline}
\newcommand{\unde}{\underline}
\newcommand{\ptf}{\pitchfork}

\newcommand{\emp}{\emptyset}
\newcommand{\wh}{\widehat}

\newcommand{\sub}{ Subsection}

\newcommand{\lc}{\lceil}
\newcommand{\rc}{\rceil}
\newcommand{\sps}{\supset}

\newcommand{\sm}{\setminus}
\newcommand{\ems}{\emptyset}
\newcommand{\sbs}{\subset}

\newcommand{\subs}{\subsection}
\newcommand{\ity}{\infty}

%%%%%%%%%%%%%%%%%%%
%%%%%%%%%%%%%%%%%%%%
%MANUSCRIPT LETTERS
%%%%%%%%%%%%%%%%%%%
%%%%%%%%%%%%%%%%%%%%

\newcommand{\GC}{\GG\CC}
\newcommand{\GCT}{\GG\CC\TT}
\newcommand{\GT}{\GG\TT}

\newcommand{\GA}{\GG\AA}
\newcommand{\GRP}{\GG\RR\PP}

\newcommand{\GgC}{\GG\gC}
\newcommand{\GgCC}{\GG\gC\CC}

\newcommand{\GgCT}{\GG\gC\TT}

\newcommand{\GgCY}{\GG\gC\YY}
\newcommand{\GgCYT}{\GG\gC\YY\TT}

\newcommand{\GCCT}{\GG\gC\CC\TT}
\newcommand{\GCC}{\GG\gC\CC}

\newcommand{\GKSC}{\GG\KK\SS\gC}
\newcommand{\GKS}{\GG\KK\SS}

%%%%%%%%%%%%%%%%%%%
%%%%%%%%%%%%%%%%%%%%
%ARROWS
%%%%%%%%%%%%%%%%%%%
%%%%%%%%%%%%%%%%%%%%

\newcommand{\strr}[3]
{{#1}^{\displaystyle\twoheadrightarrow}_{[{#2},{#3}]}}

\newcommand{\str}[1]{{#1}^{\displaystyle\twoheadrightarrow}}

\newcommand{\stind}[3]
{{#1}^{\displaystyle\rightsquigarrow}_{[{#2},{#3}]}}

\newcommand{\st}[1]{\overset{\rightsquigarrow}{#1}}
\newcommand{\bst}[1]{\overset{\displaystyle\rightsquigarrow}
\to{\boldkey{#1}}}

\newcommand{\stexp}[1]{{#1}^{\rightsquigarrow}}
\newcommand{\bstexp}[1]{{#1}^{\displaystyle\rightsquigarrow}}

\newcommand{\bstind}[3]{{\boldkey{#1}}^{\displaystyle\rightsquigarrow}_
{[{#2},{#3}]}}
\newcommand{\bminstind}[3]{\stind{({\boldkey{-}\boldkey{#1}})}{#2}{#3}}

\newcommand{\ST}{\stexp}

\newcommand{\stv}{\stexp {(-v)}}
\newcommand{\stu}{\stexp {(-u)}}
\newcommand{\stw}{\stexp {(-w)}}

\newcommand{\strv}[2]{\stind {(-v)}{#1}{#2}}
\newcommand{\strw}[2]{\stind {(-w)}{#1}{#2}}
\newcommand{\stru}[2]{\stind {(-u)}{#1}{#2}}

\newcommand{\stvv}[2]{\stind {v}{#1}{#2}}
\newcommand{\stuu}[2]{\stind {u}{#1}{#2}}
\newcommand{\stww}[2]{\stind {w}{#1}{#2}}

\newcommand{\vovo}{\stexp {(-v0)}}
\newcommand{\vov}{\stexp {(-v1)}}

\newcommand{\fl}[1]{{#1}\!\da}
\newcommand{\fll}[1]{({#1}\!\da)}

\newcommand{\vflesh}{\fl{v}}
\newcommand{\wflesh}{\fl{w}}

\newcommand{\vfllesh}{\fll{v}}
\newcommand{\wfllesh}{\fll{w}}

\newcommand{\vvfl}{\wi v\!\da}
\newcommand{\wivflesh}{\wi v\!\da}

\newcommand{\RA}{\Rightarrow}
\newcommand{\LA}{\Leftarrow}
\newcommand{\RLA}{\Leftrightarrow}
\newcommand{\LRA}{\Leftrightarrow}

\newcommand{\lau}[1]{{\xleftarrow{#1}}}

\newcommand{\rau}[1]{{\xrightarrow{#1}}}
%newcommand{\lad}[1]{{\xleftarrow[#1]{}}}
\newcommand{\rad}[1]{ {\xrightarrow[#1]{}} }

\newcommand{\da}{\downarrow}

%%%%%%%%%%%%%%%%%%%%%%%%%%%%%%%%%%%%%%%%%%%%%%%%%%%%%%%%%%%%
%%%%%%%%%%%%%%%%%%%%%%%%%%%%%%%%%%%%%%%%%%%%%%%%%%%%%%%%%%%%
%%%%%SPACES OF VECTOR FIELDS
%%%%%%%%%%%%%%%%%%%%%%%%%%%%%%%%%%%%%%%%%%%%%%%%%%%%%%%%%%%%
%%%%%%%%%%%%%%%%%%%%%%%%%%%%%%%%%%%%%%%%%%%%%%%%%%%%%%%%%%%%%

\newcommand{\vecm}{\Vect^1_0(M)}
\newcommand{\vecbn}{\Vect^1_b(N)}
\newcommand{\vecw}{\Vect^1(W)}
\newcommand{\vecrm}{\Vect^1_b(\RRR^m)}
\newcommand{\vecbm}{\Vect^1_b(M)}
\newcommand{\ver}{\text{\rm Vect}^1(\RRR^ n)}
\newcommand{\verr}{\text{\rm Vect}^1_0(\RRR^ n)}
\newcommand{\hrrr}{\text{\rm Vect}^1(M)}
\newcommand{\vemm}{\text{\rm Vect}^1_0(M)}

\newcommand{\vem}{\text{{\rm Vectt}}(M)}
\newcommand{\vebK}{\text{{\rm Vectt}}(B,K)}
\newcommand{\vemK}{\text{{\rm Vectt}}(M,K)}
\newcommand{\vemc}{\text{{\rm Vectt}}_c(M)}
\newcommand{\vemQ}{\text{{\rm Vectt}}(M,Q)}

\newcommand{\vectt}[1]{\text{{\rm Vectt}}(#1)}

\newcommand{\vew}{\text{\rm Vect}^1 (W,\bot)}

\newcommand{\veww}{\text{\rm Vect}^1 (W)}

%%%%%%%%%%%%%%%%%%%%%%%%%%%%%%%%%%%%%%%%%%%%%%%%%%%%%%%%%%%%
%%%%%%%%%%%%%%%%%%%%%%%%%%%%%%%%%%%%%%%%%%%%%%%%%%%%%%%%%%%%
%%%%%RINGS MODULES TENSOR K1
%%%%%%%%%%%%%%%%%%%%%%%%%%%%%%%%%%%%%%%%%%%%%%%%%%%%%%%%%%%%
%%%%%%%%%%%%%%%%%%%%%%%%%%%%%%%%%%%%%%%%%%%%%%%%%%%%%%%%%%%%%

\newcommand{\tens}[1]{\underset{#1}{\otimes}}

\newcommand{\LLxi}{\wi \L_\xi}

\newcommand{\Lx}{\L_{(\xi)}}
\newcommand{\Lxi}{{\wh \L}_\xi}
\newcommand{\lL}{\wh{\wh L}}

\newcommand{\Rxi}{{\ove R}_\xi}
\newcommand{\Nxi}{{\ove N}_\xi}
\newcommand{\Rcxi}{{\bar R}_\xi^c}

\newcommand{\sil}{ S^{-1}\L }
\newcommand{\kil}{\ove{K}_1(\L)}
\newcommand{\killl}{\ove{K}_1(\wh\L)}
\newcommand{\kisl}{\ove{K}_1(S^{-1}\L )}

\newcommand{\klxi}{K_1(\Lxi)}
\newcommand{\kklxi}{\ove{K_1}(\Lxi)}

\newcommand{\popo}{\tens{\L}\Lxi}
\newcommand{\popom}{\tens{\L^-}\Lxi^-}

%%%%%%%%%%%%%%%%%%%%%%%%%%%%%%%%%%%%%%%%%%%%%%%%%%%%%%%%%%%%
%%%%%%%%%%%%%%%%%%%%%%%%%%%%%%%%%%%%%%%%%%%%%%%%%%%%%%%%%%%%
%%%%%KRAYA KOBORDISMA OBOZNACH V_a, V_b A TAKZHE SAMO V
%%%%%%%%%%%%%%%%%%%%%%%%%%%%%%%%%%%%%%%%%%%%%%%%%%%%%%%%%%%%
%%%%%%%%%%%%%%%%%%%%%%%%%%%%%%%%%%%%%%%%%%%%%%%%%%%%%%%%%%%%%

\newcommand{\amk}{\AA^{(m)}_k}
\newcommand{\amkm}{\AA^{(m)}_{k-1}}
\newcommand{\bmk}{\BB^{(m)}_k}
\newcommand{\bmkm}{\BB^{(m)}_k}

\newcommand{\tivm}{\wi V^-}

\newcommand{\vk}{V_{\langle k\rangle}^-}
\newcommand{\tivkm}{\wi V_{\langle k-1\rangle}^-}
\newcommand{\tivk}{\wi V_{\langle k\rangle}^-}

\newcommand{\vkm}{V_{\langle k-1\rangle}^-}
\newcommand{\vkp}{V_{\langle k+1\rangle}^-}

\newcommand{\hvk}{\wh V_{\langle k\rangle}^-}

\newcommand{\hvkm}{\wh V_{\langle k-1\rangle}^-}
\newcommand{\hvkp}{\wh V_{\langle k+1\rangle}^-}

\newcommand{\vkvk}{V_{\prec k\succ}^-}

\newcommand{\vkvkm}{V_{\prec k-1\\succ}^-}

\newcommand{\tivkvk}{\wi V_{\prec k\succ}^-}

\newcommand{\tivkvkm}{\wi V_{\prec k-1 \succ}^-}

\newcommand{\vvvbs}{V_b^{( s)    }}
\newcommand{\vvvas}{V_a^{( s)    }}
\newcommand{\vvvbsm}{V_b^{( s-1)    }}
\newcommand{\vvvasm}{V_a^{( s-1)    }}

\newcommand{\vvbsm}{V_b^{( s-1)    }}
\newcommand{\vvasm}{V_a^{( s-1)    }}
\newcommand{\vvbs}{V_b^{[s]}    }
\newcommand{\vvas}{V_a^{[ s]    }}

\newcommand{\factor}{\vvbs / \vvbsm}
\newcommand{\factora}{\vvas / \vvasm}

\newcommand{\vvksm}{\wi V_k^{( s-1)   }}
\newcommand{\vvks}{\wi V_k^{[s]}    }

\newcommand{\vvkmsm}{\wi V_{k-1}^{( s-1)   }}
\newcommand{\vvkms}{\wi V_{k-1}^{[s]}    }

\newcommand{\fac}{\vvks / \vvksm}
\newcommand{\facm}{\vvkms / \vvkmsm}

\newcommand{\vbsm}{V_b^{\{\leq s-1\}}    }
\newcommand{\vasm}{V_a^{\{\leq s-1\}}    }
\newcommand{\vbs}{V_b^{\{\leq s\}}    }
\newcommand{\vas}{V_a^{\{\leq s\}}    }

\newcommand{\wivksm}{\wi V_k^{\{\leq s-1\}}    }
\newcommand{\wivkmsm}{\wi V_{k-1}^{\{\leq s-1\}}    }
\newcommand{\wivks}{\wi V_k^{\{\leq s\}}    }
\newcommand{\wivkms}{\wi V_{k-1}^{\{\leq s\}}    }

\newcommand{\Vbsm}{V_b^{[\leq s-1]}(\d)    }
\newcommand{\Vasm}{V_a^{[\leq s-1]}(\d)    }
\newcommand{\Vbs}{V_b^{[\leq s]}(\d)    }
\newcommand{\Vas}{V_a^{[\leq s]}(\d)    }

\newcommand{\vass}{V_{a_{s+1}}}

\newcommand{\vbkm}{V_b^{\{\leq k-1\}}    }
\newcommand{\vakm}{V_a^{\{\leq k-1\}}    }
\newcommand{\vbk}{V_b^{\{\leq k\}}    }
\newcommand{\vak}{V_a^{\{\leq k\}}    }

\newcommand{\Vbkm}{V_b^{[\leq k-1]}(\d)    }
\newcommand{\Vakm}{V_a^{[\leq k-1]}(\d)    }
\newcommand{\Vbk}{V_b^{[\leq k]}(\d)    }
\newcommand{\Vak}{V_a^{[\leq s]}(\d)    }

%%%%%%%%%%%%%%%%%%%%%%%%%%%%%%%%%%%%%%%%%%%%%%%%%%%%
%%%%%%%%%%%%%%%%%%%%%%%%%%%%%%%%%%%%%%%%%%%%%%%%%%%%
%%%%%KRAYA KOBORDISMA OBOZNACH \PARTIAL
%%%%%%%%%%%%%%%%%%%%%%%%%%%%%%%%%%%%%%%%%%%%%%%%%%%%
%%%%%%%%%%%%%%%%%%%%%%%%%%%%%%%%%%%%%%%%%%%%%%%%%%%%

\newcommand{\dow}{\pr_0 W}

\newcommand{\daw}{\pr_1 W}

\newcommand{\hdaw}{\wh{\pr_1 W}}

\newcommand{\tipwk}{(\pr_1 \wi W)^{\{\leq k\}}}

\newcommand{\timwk}{(\pr_0 \wi W)^{\{\leq k\}}}

\newcommand{\tipwkm}{(\pr_1 \wi W)^{\{\leq k-1\}}}

\newcommand{\pws}{(\pr_1 W)^{\{\leq s\}}}

\newcommand{\hpws}{(\wh{\pr_1 W})^{\{\leq s\}}}

\newcommand{\hpwsm}{(\wh{\pr_1 W})^{\{\leq {s-1}\}}}

\newcommand{\mws}{(\pr_0 W)^{\{\leq s\}}}

\newcommand{\pwk}{(\pr_1 W)^{\{\leq k\}}}

\newcommand{\mwk}{(\pr_0 W)^{\{\leq k\}}}

\newcommand{\pwkm}{(\pr_1 W)^{\{\leq k-1\}}}

\newcommand{\pwsm}{(\pr_1 W)^{\{\leq s-1\}}}

\newcommand{\pwkmm}{(\pr_1 W)^{\{\leq k-2\}}}

\newcommand{\mwkmm}{(\pr_0 W)^{\{\leq k-2\}}}

\newcommand{\mwkm}{(\pr_0 W)^{\{\leq k-1\}}}

\newcommand{\mwsm}{(\pr_0 W)^{\{\leq s-1\}}}

\newcommand{\mwkp}{(\pr_0 W)^{\{\leq k+1\}}}

\newcommand{\dwmok}{\daw^{( k)}}

\newcommand{\dwmokp}{\daw^{(k+1)}}

\newcommand{\dwmokm}{\daw^{(k-1)}}
\newcommand{\dwmokmm}{\daw^{(k-2)}}

%%%%%%%%%%%%%%%%%%%%%%%%%%%%%%%%%%%%%%%%%%%%%%%%%%%%
%%%%%%%%%%%%%%%%%%%%%%%%%%%%%%%%%%%%%%%%%%%%%%%%%%%%
%%%%%SAM  KOBORDISM W I EGO FILTRATSII
%%%%%%%%%%%%%%%%%%%%%%%%%%%%%%%%%%%%%%%%%%%%%%%%%%%%
%%%%%%%%%%%%%%%%%%%%%%%%%%%%%%%%%%%%%%%%%%%%%%%%%%%%

\newcommand{\moi}[1]{\MM^{(0)}_{#1}}
\newcommand{\moii}[1]{\MM^{(1)}_{#1}}

\newcommand{\Wal}{W_{[a,\l]}}

\newcommand{\Wlm}{W_{[\l,\m]}}
\newcommand{\Wlb}{W_{[\l,b]}}

\newcommand{\Wam}{W_{[a,\m]}}

\newcommand{\Wall}{W_{[a,\l']}}

\newcommand{\Wkr}{W^{\circ}}
\newcommand{\wkr}{W^{\circ}}

\newcommand{\wa}[2]{ W_{[a_{#1}, a_{#2}]}}

\newcommand{\waa}[1]{ W_{[a, a_{#1}]}}

\newcommand{\Wa}[2]{ W_{[{#1}, {#2}]}}

\newcommand{\WS}[1]{ W^{\{\leq {#1}\}}}

\newcommand{\ws}{\WS {s}}

\newcommand{\wsm}{\WS {s-1}}

\newcommand{\wsmm}{\WS {s-2}}

\newcommand{\wk}{\WS {k}}

\newcommand{\wkm}{\WS {k-1}}

\newcommand{\wkmm}{\WS {k-2}}

\newcommand{\wsn}{ W^{[\leq s]}(\nu)}

\newcommand{\wsmn}{ W^{[\leq s-1]}(\nu)}

\newcommand{\wsk}{ W^{[\leq k]}(\nu)}

\newcommand{\Wmok}{W^{\prec k\succ}}
\newcommand{\wmok}{W^{\langle k\rangle}}

\newcommand{\wmokp}{W^{\langle k+1\rangle}}

\newcommand{\wmokm}{W^{\langle k-1\rangle}}
\newcommand{\wmokmm}{W^{\langle k-2\rangle}}

\newcommand{\wmos}{W^{\langle s\rangle}}

\newcommand{\wmosm}{W^{\langle s-1\rangle}}

\newcommand{\wmoo}{W^{\langle 0\rangle}}

\newcommand{\wwk}{\( \wmok , \wmokm \)}

\newcommand{\wwkp}{\( \wmokp , \wmok \)}

\newcommand{\wwkm}{\( \wmokm , \wmokmm \)}

\newcommand{\wws}{\bigg( \wmos , \wmosm \bigg)}

\newcommand{\wasn}{W^{\lc s\rc}(\nu)}

\newcommand{\wakn}{W^{\lc k\rc}(\nu)}

\newcommand{\hwm}{H_*\( \wmok, \wmokm\)}

\newcommand{\hkwm}{H_k\( \wmok, \wmokm\)}

\newcommand{\hwmp}{H_*\( \wmokp, \wmok\)}

\newcommand{\hkwmm}{H_k\( \wmokm, \wmokmm\)}

\newcommand{\hkwmp}{H_k\( \wmokp, \wmok\)}

\newcommand{\tiws}{\wi W^{\{\leq s\}}}

\newcommand{\tiwk}{\wi W^{\{\leq k\}}}

\newcommand{\tiwsm}{\wi W^{\{\leq s-1\}}}

\newcommand{\tiwkm}{\wi W^{\{\leq k-1\}}}

\newcommand{\wkwk}{ W^{[ k]}}
\newcommand{\tiwmok}{\wi W^{\prec k\succ}}
\newcommand{\womok}{W_0^{\prec k\succ}}

\newcommand{\Womok}{W_0^{\langle k\rangle}}

\newcommand{\Wmokp}{W^{\langle k+1\rangle}}

\newcommand{\Wmokm}{W^{\langle k-1\rangle}}
\newcommand{\Wmokmm}{W^{\langle k-2\rangle}}

\newcommand{\Wmos}{W^{\langle s\rangle}}

\newcommand{\Wmosm}{W^{\langle s-1\rangle}}

\newcommand{\Wmoo}{W^{\langle 0\rangle}}

\newcommand{\hWmok}{\wh W^{\langle k\rangle}}
\newcommand{\hWmokm}{\wh W^{\langle k-1\rangle}}

\newcommand{\tiWmok}{\wi W^{\langle k\rangle}}
\newcommand{\tiWmokm}{\wi W^{\langle k-1\rangle}}

%%%%%%%%%%%%%%%%%%%%%%%%%%%%%%%%%%%%%%%%%%%%%%%%%%%%
%%%%%%%%%%%%%%%%%%%%%%%%%%%%%%%%%%%%%%%%%%%%%%%%%%%%
%%%%%%%%%%%%%%%%%%%%%%%%%%%%%%%%%%%%%%%%%%%%%%%%%%%%

%%%%%%%%%%%%%%%%%%%%%%%%%%%%%%%%%%%%%%%%%%%%%
%%%%%%%%%%%%%%%%%%%%%%%%%%%%%%%%%%%%%%%%%%%%%%
%MACROS_SLOVA
%%%%%%%%%%%%%%%%%%%%%%%%%%%%%%%%%%%%%%%%%%%%%
%%%%%%%%%%%%%%%%%%%%%%%%%%%%%%%%%%%%%%%%%%%%%%

\newcommand{\ifff}{if and only if}

\newcommand{\orial}{oriented almost transverse}
\renewcommand{\th}{therefore}
\newcommand{\at}{almost~ transverse}
\newcommand{\ata}{almost~ transversality~ condition}
\newcommand{\gr}{gradient}
\newcommand{\Mf}{Morse function}
\newcommand{\iis}{it is sufficient}
\newcommand{\sut}{~such~that~}
\newcommand{\sufsm}{~sufficiently~ small}
\newcommand{\sufla}{~sufficiently~ large}
\newcommand{\sufcl}{~sufficiently~ close}
\newcommand{\wrt}{~with respect to}
\newcommand{\ho}{homomorphism}
\newcommand{\iso}{isomorphism}
\newcommand{\rgr}{Riemannian gradient}
\newcommand{\riemm}{Riemannian metric}

\newcommand{\trasp}{trajectory starting at a point of}
\newcommand{\trasps}{trajectories starting at a point of}

\newcommand{\ma}{manifold}
\newcommand{\nei}{neighborhood}
\newcommand{\dfm}{diffeomorphism}

\newcommand{\vf}{vector field}

\newcommand{\vfs}{vector fields}

\newcommand{\fe}{for every}

\renewcommand{\top}{topology}

\newcommand{\tr}{~trajectory }

\newcommand{\grs}{~gradients}
\newcommand{\trs}{~trajectories}

\newcommand{ \co}{~cobordism}
\newcommand{
\sma}{submanifold}
\newcommand{
\hos}{~homomorphisms}
\newcommand{
\Th}{~Therefore}

\newcommand{
\tthen}{\text \rm ~then}

\newcommand{
\wwe}{\text \rm ~we  }
\newcommand{
\hhave}{\text \rm ~have}
\newcommand{
\eevery}{\text \rm ~every}

\newcommand{\noconf}{~there~is~no~possibility~of~confusion}

\newcommand{\ATA}{almost~ transversality~ tondition}
\newcommand{\cob}{~cobordism}

\newcommand{\hot}{homotopy}

\newcommand{\emem}{elementary modification}
\newcommand{\emems}{elementary modifications}

\newcommand{\TA}{transversality condition}

\newcommand{\hog}{homology}

\newcommand{\cog}{cohomology}

\newcommand{\wat}{ We shall assume that}

\newcommand{\sclv}{sufficiently close to $v$ in $C^0$-topology}

\newcommand{\cf}{continuous function }

\newcommand{\heq}{homotopy equivalence}

\newcommand{\heeq}{homology equivalence}

\newcommand{\eg}{exponential growth}

\newcommand{\nics}{Novikov incidence coefficients}
\newcommand{\nic}{Novikov incidence coefficient}

\newcommand{\negc}{Novikov exponential growth conjecture}

\newcommand{\mc}{Morse Complex   }

\newcommand{\mas}{manifolds   }

\newcommand{\nc}{Novikov Complex   }

\newcommand{\glvf}{gradient-like vector field}

\newcommand{\glvfs}{gradient-like vector fields}

\newcommand{\fg}{finitely generated   }

\newcommand{\she}{simple~homotopy~equivalence}

\newcommand{\sht}{simple~homotopy~type}

\newcommand{\ta}{transversality condition}

\newcommand{\cpc}{convex polyhedral cone}
\newcommand{\rcpc}{rational convex polyhedral cone}

\newcommand{\mnp}{Morse-Novikov pair}

\newcommand{\rp}{rationality property}
\newcommand{\wvf}{Whitney vector field}

\newcommand{\egp}{exponential growth property}

\newcommand{\lzf}{Lefschetz zeta function}

\newcommand{\babs}{by abuse of notation}
\newcommand{\su}{subsection}
\newcommand{\Prop}{\text{Proposition}}

\newcommand{\aand}{\quad\text{and}\quad}
\newcommand{\wwhere}{\quad\text{where}\quad}
\newcommand{\ffor}{\quad\text{for}\quad}
\newcommand{\iif}{\quad\text{if}\quad}
\newcommand{\iiif}{~\text{if}~}

\newcommand{\eqi}{equivalence}

\newcommand{\mfcobv}{~ Let $\fcob$
be a Morse function on
 a cobordism $W$ and $v$
 be an $f$-gradient. ~}

\newcommand{\mfcob}{~ Let $\fcob$
be a Morse function on
 a cobordism $W$}
\newcommand{\msf}{Morse-Smale filtration}

\newcommand{\fbfg}{ free based finitely generated }
%%%%%%%%%%%%%%%%%%%%%%%%%%%%%%%%%%%%%%%%%%%%%%%%%%%%%%
%%%%%%%%%%%%%%%%%%%%%%%%%%%%%%%%%%%%%%%%%%%%%%%%%%%%%%
%%%%%%%RAZNOE
%%%%%%%%%%%%%%%%%%%%%%%%%%%%%%%%%%%%%%%%%%%%%%%%%%%%%%
%%%%%%%%%%%%%%%%%%%%%%%%%%%%%%%%%%%%%%%%%%%%%%%%%%%%%%

\newcommand{\tret}{{\frac 13}}
\newcommand{\dvet}{{\frac 23}}
\newcommand{\polt}{{\frac 32}}
\newcommand{\polo}{{\frac 12}}

\newcommand{\bv}{B(v,\d)}

\newcommand{\ti}{\times}

\newcommand{\FR}{{\mathcal{F}}r}
\newcommand{\gt}{{\mathcal{G}}t}

\newcommand{\en}{enumerate}

\newcommand{\Prf}{{\it Proof.\quad}}
\newcommand{\prf}{{\it Proof:\quad}}

\newcommand{\nr}{\Vert}
\newcommand{\smo}{C^{\infty}}

\newcommand{\fpr}[2]{{#1}^{-1}({#2})}
\newcommand{\sdvg}[3]{\widehat{#1}_{[{#2},{#3}]}}
\newcommand{\disc}[3]{B^{({#1})}_{#2}({#3})}
\newcommand{\Disc}[3]{D^{({#1})}_{#2}({#3})}
\newcommand{\desc}[3]{B_{#1}(\leq{#2},{#3})}
\newcommand{\Desc}[3]{D_{#1}(\leq{#2},{#3})}
\newcommand{\komp}[3]{{\bold K}({#1})^{({#2})}({#3})}
\newcommand{\Komp}[3]{\big({\bold K}({#1})\big)^{({#2})}({#3})}

\newcommand{\ran}{\{(A_\lambda , B_\lambda)\}_{\lambda\in\Lambda}}
\newcommand{\rran}{\{(A_\lambda^{(s)},
 B_\lambda^{(s)}  )\}_{\lambda\in\Lambda, 0\leq s\leq n }}
\newcommand{\brs}{\rran}
\newcommand{\rans}{\{(A_\sigma , B_\sigma)\}_{\sigma\in\Sigma}}

\newcommand{\fmin}{F^{-1}}
\newcommand{\vh}{\widehat{(-v)}}

\newcommand{\chart}{\Phi_p:U_p\to B^n(0,r_p)}
\newcommand{\atlas}{\{\Phi_p:U_p\to B^n(0,r_p)\}_{p\in S(f)}}
\newcommand{\flow}{{\VV}=(f,v, \UU)}

\newcommand{\Rn}{\bold R^n}
\newcommand{\Rk}{\bold R^k}

\newcommand{\fcob}{f:W\to[a,b]}

\newcommand{\phicob}{\phi:W\to[a,b]}

\newcommand{\crr}{p\in S(f)}
\newcommand{\nrv}{\Vert v \Vert}
\newcommand{\nrw}{\Vert w \Vert}
\newcommand{\nru}{\Vert u \Vert}

\newcommand{\obb}{\cup_{p\in S(f)} U_p}
\newcommand{\proob}{\Phi_p^{-1}(B^n(0,}

\newcommand{\indl}[1]{{\scriptstyle{\text{\rm ind}\leqslant {#1}~}}}
\newcommand{\inde}[1]{{\scriptstyle{\text{\rm ind}      =   {#1}~}}}
\newcommand{\indg}[1]{{\scriptstyle{\text{\rm ind}\geqslant {#1}~}}}

\newcommand{\obbi}{\cup_{p\in S_i(f)}}

\newcommand{\pr}{\partial}
\newcommand{\prx}[2]{\frac {\pr {#1}}{\pr x} ({#2})}
\newcommand{\pry}[2]{\frac {\pr {#1}}{\pr y} ({#2})}
\newcommand{\prz}[2]{\frac {\pr {#1}}{\pr z} ({#2})}
\newcommand{\przbar}[2]{\frac {\pr {#1}}{\pr \bar z} ({#2})}
\newcommand{\chape}[2]{\frac  {\pr {#1}}{\pr {#2}} }
\newcommand{\chapee}[2]{\frac  {\pr^2 {#1}}{\pr {#2}^2} }

\newcommand{\xit}{\tilde\xi_t}

\newcommand{\VODIN}{V_{1/3}}
\newcommand{\VDVA}{V_{2/3}}
\newcommand{\VM}{V_{1/2}}
\newcommand{\ddd}{\cup_{p\in S_i(F_1)} D_p(u)}
\newcommand{\dddmin}{\cup_{p\in S_i(F_1)} D_p(-u)}
\newcommand{\where}{\quad\text{\rm where}\quad}

\newcommand{\kr}[1]{{#1}^{\circ}}

\newcommand{\mods}{\vert s(t)\vert}
\newcommand{\exd}{e^{2(D+\alpha)t}}
\newcommand{\exmin}{e^{-2(D+\alpha)t}}

\newcommand{\intt}{[-\theta,\theta]}

\newcommand{\ffmin}{f^{-1}}

\newcommand{\vxi}{v\langle\vec\xi\rangle}

\newcommand{\qt}{\hfill\triangle}
\newcommand{\qs}{\hfill\square}

\newcommand{\pa}{\vskip0.1in}

\renewcommand{\(}{\big(}
\renewcommand{\)}{\big)}

\newcommand{\Vm}{V_\m}

\newcommand{\Vl}{V_\l}

\newcommand{\lccc}{\wh\L_{C}}

\newcommand{\ld}{\wh\L_{D}}

\newcommand{\udp}{{\displaystyle {\vartriangle}}}
\newcommand{\ddp}{{\displaystyle {\triangledown}}}

\newcommand{\Vv}{{\boldsymbol{v}}}

\newcommand{\hV}{\wh V}
\newcommand{\hHH}{\wh \HH}

\newcommand{\gama}[2]{\g({#1}, \tau_a({#2},{#1}); w )}

\newcommand{\gam}[2]{\g({#1}, \tau_0({#2},{#1}); w )}
\newcommand{\ga}[2]{\g({#1}, \tau({#2},{#1}); w )}

\newcommand{\mi}[3]{{#1}^{-1}\([{#2},{#3}]\)}

\newcommand{\fii}[2]{\mi {\phi}{a_{#1}}{a_{#2}} }

\newcommand{\fifi}[2]{\mi {\phi}{#1}{#2} }

\newcommand{\pf}[2]{\mi {\phi_1}{\a_{#1}}{\a_{#2}} }

\newcommand{\mf}[2]{\mi {\phi_0}{\b_{#1}}{\b_{#2}}}

\newcommand{\dqr}{\pr_- Q_r}

\newcommand{\ds}{\pr_s}

\newcommand{\dsm}{\pr_{s-1}}

\newcommand{\yz}{Y_k(v)\cup Z_k(v)}

\newcommand{\Gama}{{\nazad{ \Gamma}}}
\newcommand{\ug}[1]{\llcorner {#1} \lrcorner}
\newcommand{\npqv}{n(\bar p, \bar q; v)}
\newcommand{\fms}{f:M\to S^1   }

\newcommand{\nkpqv}{n_k(\bar p, \bar q; v)}

\newcommand{\GLT}{\GG lt}

\newcommand{\Trln}{{\text Trln}}

\newcommand{\Trlln}{{\text TrLn}}

\newcommand{\Tr}{{\text{\rm  Tr}}}
\newcommand{\TrL}{{\text TrL}}
\newcommand{\limdir}{\underset {\to}{\lim}}
\newcommand{\liminv}{\underset {\leftarrow}{\lim}}

\newcommand{\kom}[2]{ {#1}{#2}{ {#1}^{-1}} {{#2}^{-1}} }

\newcommand{\komm}[2]{ {#1}{#2}{ ({#1})^{-1}} {({#2})^{-1}} }
\newcommand{\kommm}[2]{ {#1}'{#2}'{ ({#1}'')^{-1}} {({#2}'')^{-1}} }

\newcommand{\Trll}{\TL'}

\newcommand{\cmd}{ C_*^\D( \wi M)}
\newcommand{\cmxi}{\wh C_*^\D( \wi M, \xi)}
\newcommand{\whgxi}{\wh {{\rm Wh}} (G,\xi)}

\newcommand{\vwdwp}{Vect(W,\pr_0W;P)}
\newcommand{\ewdwp}{\EE(W,\pr_0W)}
\newcommand{\ewdwo}{\EE(W_1,\pr_0W_1)}
\newcommand{\ewdwd}{\EE(W_2,\pr_0W_2)}

\newcommand{\kpr}{K_r^+}

\newcommand{\kmr}{K_r^-}

\newcommand{\kpd}{K_r^+(\d)}

\newcommand{\kmd}{K_r^-(\d)}

\newcommand{\addc}{\addtocontents{toc}{\protect\vspace{10pt}}}

\newcommand{\mxi}{M_\xi   }

\newcommand{\cmm}{C_*^\D(\wi M)}

\newcommand{\cvm}{C_*^\D(\wi V^-)}

\newcommand{\ey}{\wi E_*}
\newcommand{\eey}{\wi \EE_*}

\newcommand{\eky}{\wi E(k)_*}

\newcommand{\eti}{\wi{\wi\EE_*}}

\newcommand{\etik}{\wi{\wi\EE}_k}

\newcommand{\etikp}{\wi{\wi\EE}_{k+1}}

\newcommand{\ctiu}{\wi{\wi C}_*(u_1)}

\newcommand{\ctiuk}{\wi{\wi C}_k(u_1)}

\newcommand{\ctiv}{\wi{\wi C}_*(v)}

\newcommand{\ctiukm}{\wi{\wi C}_{k-1}(u_1)}

\newcommand{\scc}[1]{|{\scriptscriptstyle{#1}}}\newcommand{\rrr}{\{\wi r\}}

\newcommand{\tidow}{\pr_0 \wi W}
\newcommand{\tidaw}{\pr_1 \wi W}

\newcommand{\tivkp}{\wi V_{\langle k+1\rangle}^-}

\newcommand{\ur}[1]{\overset{\smallfrown}{#1}}

\newcommand{\dr}[1]{\underset{\smallsmile}{#1}}

\newcommand{\uUu}{\overset{\twoheadrightarrow}{u}}
\newcommand{\vVv}{\overset{\twoheadrightarrow}{v}}
\newcommand{\wWw}{\overset{\twoheadrightarrow}{w}}

\newcommand{\bfun}{{\bf 1}}

\newcommand{\ppmm}{{\scriptstyle{ \pm}}}

%%%%%%%%%%%%%%%%%%%%%%%%%%%%%%%%%%%%%%%%%%%%%
%%%%%%%%%%iz n_a
%%%%%%%%%%%%%%%%%%%%%%%%%%

%\newcommand{\Lxi}{{\wh \L}_\xi}
\newcommand{\Lxim}{{\wh \L}^-_\xi}

\newcommand{\lxi}{{\bar \L}_\xi}

\newcommand{\Lc}{{\wh \L}_C}
\newcommand{\Lcm}{{\wh \L}_C^-}

\newcommand{\lLL}{\wh{\wh \L}}

\newcommand{\Xc}{{\wh X}_C}
\newcommand{\Xfaa}{{\wh X}_{(F_i,\vec\a)}}

\newcommand{\bs}{\boldsymbol}

%%%%%%%%%%%%%%%%%%
%%%%%%%skobka
%%%%%%%%%%%%%%

\newcommand{\bikl}{\text{\rM (}}
\newcommand{\bikr}{\text{\rM )}}

\newcommand{\ck}{C_*^{(k)}}
\newcommand{\dk}{D_*^{(k)}}
\newcommand{\ckp}{C_*^{(k+1)}}
\newcommand{\dkp}{D_*^{(k+1)}}

%ii\input /home/a/bibE
%%%%%%%%%%%%%%

%%%%%%%%%%%%%%%%%%%%%%%%%%%%%%%%%%%%
%ARTICLES IN FOREIGN LANGUAGES
%%%%%%%%%%%%%%%%%%%%%%%%%%%%%%%%%%%%

\newcommand{\arnoldXC}{ V.I. Arnold,
\emph{Dynamics of intersections},
        Proceedings of a Conference
 in Honour of
J.Moser, edited by
 P.Rabinowitz and R.Hill,
 Academic Press,
 1990
pp. 77--84.  }

\newcommand{\arnoldXCprim}{ V.I. Arnold,
\emph{Dynamics of Complexity of Intersections},
Boletim SOc. Brasil. Mat. (N.S.),
1990, {\bff 21} (1), 1-10.}

\newcommand{\arnoldXCIII}{ V.I. Arnold,
\emph{Bounds for Milnor numbers of intersections
in holomorphic dynamical systems}, In:
Topological Methods in Modern Mathematics,
 Publish or Perish,
 1993,
pp. 379--390.}

\newcommand{\arnoldXCIV}{ V.I. Arnold,
\emph{Sur quelques probl\`emes de
 la th\'eorie des syst\`emes
dynamiques},
  Journal of the Julius Schauder center, {\bff 4}
 1994
pp. 209--225.  }

\newcommand{\armaz}
{M.Artin, B.Mazur,
\emph{On periodic points},
Annals of Math. {\bff 102} (1965), 82--99.
}

\newcommand{\baladi}
{V.Baladi,
\emph{Periodic orbits and dynamical spectra},
Ergodic theory and dynamical spectra,
{\bff 18}
(1998),
    255 - 292.
}

\newcommand{\chapman}
{ T.A.Chapman,
\emph{ Topological invariance of Whitehead torsion},
American J. of Math.
{\bff 96}
(1974),
    488 - 497
}

\newcommand{\bhs}
{H.Bass, A.Heller, R.G.Swan,
\emph{The Whitehead group of a polynomial extension},
Inst. Hautes Etudes Sci. Publ. Math. {\bff 22} (1964),
61--79
}

\newcommand{\brlev}
{ W.Browder, J.Levine,
\emph{Fibering manifolds over the circle},
Comment. Math. Helv. {\bff 40} (1966),
153--160
}

\newcommand{\farran}{ M. Farber and A. A. Ranicki,
\emph{ The Morse-Novikov theory of circle-valued functions
and noncommutative localization,}

  Proc. 1998 Moscow Conference
for the 60th Birthday of S. P. Novikov, tr.
 Mat. Inst. Steklova, {\bf 225}, 1999,
381 -- 388.

\quad E-print:

dg-ga/9812122.
}

\newcommand{\farrell}
{F.T.Farrell,
\emph{The obstruction to fibering a manifold over a circle},
Indiana Univ.~J. {\bff 21} (1971), 315--346.
}

\newcommand{\farhsi}
{F.T.Farrell, W.-C.Hsiang,
\emph{A formula for $K_1R_\alpha[T]$},
Proc. Symp. Pure Math., Vol. {\bff 17} (1968), 192--218}

\newcommand{\fel}{  A.Fel'shtyn,
\emph{
Dynamical zeta functions,
Nielsen Theory and Reidemeister torsion,}
preprint ESI 539 of The Erwin Schr\"odinger International Institute for
Mathematical Physics,
(to appear in Memoirs of AMS)
(1998)
}

\newcommand{\franks}
{   J.Franks
\emph{Homology and dynamical systems},
CBMS Reg. Conf. vol. 49, AMS, Providence 1982.
}

\newcommand{\fuller}{  F.B.Fuller,
\emph{An index of fixed point type for periodic orbits},
Amer. J.Math
{\bff 89},
(1967)
133--148
}

\newcommand{\fried}{  D.Fried,
\emph{Homological Identities for closed orbits},
Inv. Math. {\bff 71}, (1983) 419--442.
}

\newcommand{\friedtwi}{  D.Fried,
\emph{Periodic points and twisted coefficients},
Lect. Notes in Math.,
{\bff 1007},
(1983)
261--293.
}

\newcommand{\friednewzeta}{  D.Fried,
\emph{Flow equivalence, hyperbolic systems and a new zeta function
for flows},
Comm. Math. Helv.,
{\bff 57},
(1982)
237--259.
}

\def\gnXCIV
{  R.Geoghegan, A.Nicas,
\emph{Trace and torsion in the theory of flows},
Topology,
{\bff 33},
(1994)
683--719
}

\def\gnXCIVprim
{  R.Geoghegan, A.Nicas,
\emph{Parameterized Lefschetz-Nielsen fixed
 point theory and Hochshild homology traces},
Amer. J. Math.,
{\bff 116},
(1994)
397--446 }

\def\gnXCV
{  R.Geoghegan, A.Nicas,
\emph{Higher Euler characteristics 1},
L'Enseignement Math\'ematique,
{\bff 41},
(1995)
3--62
}

\def\hadamard
{ J.Hadamard,
\emph{Sur l'it\'eration et les solutions asymptotiques
des \'equations diff\'erentielles},
Bull. Soc. Math. France,
{\bff 29},
(1901)
224 -- 228
}

\def\higman
{ G.Higman,
\emph{Units in group rings},
Proc. London Math. Soc.,
{\bff 46},
(1940)
231 -- 248
}

\newcommand{\hulee}
{M.Hutchings, Y-J.Lee
\emph{Circle-valued Morse theory, Reidemeister torsion
and Seiberg-Witten invariants of 3-manifolds},
\quad E-print:
 dg-ga/9612004 3  Dec 1996,
 journal publication:
Topology,
{\bff 38},
(1999),
861 -- 888.
}

\newcommand{\huleee}{ M.Hutchings, Y-J.Lee
\emph{  Circle-valued Morse theory and Reidemeister torsion
},
Geometry and Topology,
{\bff 3},
(1999),
369 -- 396
}

\newcommand{\irwinP}{   M.Irwin,
\emph{ On the Stable Manifold Theorem
},
Bull. London Math. Soc,
{\bff 2},
(1970),
196 -- 198
}

\newcommand{\jiang}{  B.Jiang,
\emph{ Estimation of the number of periodic orbits}
Preprint of Universit\"at Heidelberg,
Mathematisches Institut, Heft 65, Mai 1993,
Pac. J. Math. 172, No.1, 151-185 (1996).
}

\newcommand{\kupka}
{   I.Kupka,
\emph{
Contribution \`a la th\'eorie des champs g\'en\'eriques},
Contributions to Differential equations,
{\bff 2}
    (1963    ), 457--484,
    {\bff 3}
    (1964    ), 411--420.
 }

\newcommand{\latour}
{F.Latour,
\emph{
Existence de 1-formes ferm\'ees non
singuli\`eres dans une classe de cohomologie de de Rham},
Publ. IHES {\bff 80}
(1995),
}

\newcommand{\laudenbach}
{F.Laudenbach,
\emph{On the Thom-Smale complex},
Asterisque,{\bff  205} (1992)
p. 219 -- 233.}

\newcommand{\laudsiko}
{   F.Laudenbach, J.-C.Sikorav,
\emph{
Persistance d'intersection avec la section
nulle au cours d'une isotopic hamiltonienne
 dans un fibre cotangent
},
Invent.~Math.~
{\bff 82}
    (1985     ),
pp. 349--357.
 }

\newcommand{\lueckzeta}
{ W.L\"uck,
\emph{The Universal Functorial Lefschetz Invariant }
Fundam. Math. 161, No.1-2, 167-215 (1999).

Preprint: (1998)
}

\newcommand{\milncyccov}
{J.Milnor,
\emph{ Infinite cyclic coverings},
In: Conference on the topology of manifolds,
(1968)}

\newcommand{\milnWT}
{J.Milnor,
\emph{ Whitehead Torsion},
Bull. Amer. Math. Soc.
{\bff 72}
(1966),
358 - 426.
}

\newcommand{\peixoto}
{Peixoto,
\emph{ On an approximation theorem of Kupka and Smale},
J. Diff.Eq.
{\bff 3}
(1966),
423 -- 430.
}

\newcommand{\pitcher}
{E.Pitcher,
\emph{ Inequalities of critical point theory},
Bull. Amer. Math. Soc.
{\bff 64}
(1958),
1-30.
}

\newcommand{\noviquasi}
{S.P.Novikov,
\emph{Quasiperiodic Structures in topology},
in the  book:  Topological Methods in Modern Mathematics,
 Publish or Perish,
 1993,
pp. 223--235.  }

\newcommand{\pozniakXC}
{M.Pozniak,
\emph{Triangulation of compact smooth manifolds and Morse theory}
(University of Warwick preprint, 11/1990,
published posthumously as a part of the thesis of M.Pozniak
in Translations of AMS, 2000)}

\newcommand{\pozniakXCI}
{M.Pozniak,
\emph{The Morse complex, Novikov Cohomology and Fredholm Theory}
(University of Warwick preprint, 08/1991,
published posthumously as a part
of the thesis of M.Pozniak
in Translations of AMS, 2000)}

\def\ranXCV
{ A.Ranicki,
\emph{Finite domination
and  Novikov rings},
Topology,
 {\bff 34} (1995), 619--632.}

\def\ranXCIX
{ A.Ranicki,
\emph{The algebraic construction
of the Novikov complex of a circle-valued
Morse function},

\quad  E-print: math.AT/9903090}

\newcommand{\reidemeister}
{K.Reidemeister,
\emph{Homotopieringe und Linsenr\"aume},
Hamburger Abhandl. {\bff 11} (1938), 102--109.}

\def\schuetzC
{  D.Sch\"utz,
\emph{Gradient flows of closed 1-forms and
their closed orbits},
e-print:
math.DG/0009055,
~ journal article:
Forum Math. 14(2002) 509--537.
}

\def\schuetzCI
{  D.Sch\"utz,
\emph{One-parameter fixed point theory and gradient flows
of closed 1-forms},
e-print:
math.DG/0104245,
~ journal article:
K-theory, 25(2002), 59-97.
}

\newcommand{\sieben}
{L.Siebenmann,
\emph{A total Whitehead torsion obstruction to fibering over the circle},
Comment. Math. Helv. {\bff 45} (1970), 1--48.}

\newcommand{\sikoens}
{ J.-Cl.~Sikorav,
\emph{Un probleme de disjonction par
isotopie symplectique dans un
fibr\'e cotangent},
 Ann.~Scient.~Ecole~Norm.~Sup.,{\bff 19}
 (1986),  543--552.}

\newcommand{\sikoravthese}
{ J.-Cl.~Sikorav,
\emph{
Points fixes de diff\'eomorphismes
symplectiques, intersections de sous-vari\'et\'es
lagrangiennes, et singularit\'es de un-formes ferm\'ees
}
Th\'ese de Doctorat d'Etat Es
Sciences Math\'ematiques,
Universit\'e Paris-Sud, Centre d'Orsay, 1987}

\newcommand{\smale}
{  S.~Smale, \emph{On the structure of manifolds},
 Am.~J.~Math., {\bff 84} (1962)  387--399.}

\newcommand{\smalehcob}
{  S.~Smale, \emph{Generalized Poincar\'e
onjecture in dimensions greater than four},
 Ann.~Math., {\bff 74} (1961)  391--406.}

\newcommand{\smdyn}
{  S.~Smale,
\emph{Differential dynamical systems},
 Bull. Amer. Math. Soc. {\bff 73} (1967)
747--817.}

\newcommand{\smaletransv}
{  S.~Smale,
\emph{Stable \ma s for differential
equations and diffeomorphisms
},
Ann. Scuola Norm. Superiore Pisa, {\bff 18} (1963)
97--16.}

\newcommand{\smapoi}
{  S.~Smale,
\emph{Generalized Poincare's conjecture in dimensions
greater than four},
 Ann.~Math.,
{\bff 74} (1961)
 391--406.}

\newcommand{\stal}
{J.Stallings ,
\emph{On fibering certain 3-manifolds},
 Proc. 1961 Georgia conference on the Topology of 3-manifolds,
Prentice-Hall, 1962, pp. 95--100.}

\newcommand{\thom}
{ R.Thom,
\emph{Sur une partition en cellules associ\'ee
\`a une fonction sur une vari\'et\'e},
     Comptes Rendus de l'Acad\'emie de Sciences,
{\bff 228}
(1949),
 973--975.
}

\newcommand{\witten}
{ E.Witten,
 \emph{Supersymmetry and Morse theory},
 Journal of Diff.~Geom.,
{\bff 17} (1985)
 no. 2.
}

\newcommand{\whitehead}
{ J.H.C.Whitehead
 \emph{Simple homotopy types},
 Amer. J. Math.,
{\bff 72} (1952)
 pp. 1- 57
}

%%%%%%%%%%%%%%%%%%%%%%%%%%%%%%%
%    MY ARTICLES AND PREPRINTS
%%%%%%%%%%%%%%%%%%%%%%%%%%%%%

\newcommand{\patou}
{ A.V.Pajitnov, \emph{ On the Novikov
complex for rational Morse forms},
\quad preprint:
Institut for Matematik og datalogi, Odense Universitet
Preprints 1991, No 12, Oct. 1991;
~ journal article:
Annales de la Facult\'e de Sciences de
Toulouse {\bf 4}  (1995), 297--338.
}

\newcommand{\pasur}
{ A.V.Pajitnov,
\emph{
Surgery on the Novikov Complex},
Preprint:
Rapport de Recherche CNRS URA 758,  Nantes,   1993;
journal article:
K-theory {\bff 10} (1996),  323-412.
 }

\newcommand{\pamrl}
{  A.V.Pajitnov,
\emph{Rationality and exponential growth
properties of the boundary operators in the Novikov
Complex},
Mathematical Research Letters,
{\bff 3}
(1996),
  541-548.
 }

\newcommand{\paasym}
{  A.V.Pajitnov,
\emph{   On the asymptotics of
Morse numbers of finite covers of
manifold
},

\quad E-print:
math.DG/9810136, 22 Oct 1998,\quad
journal article:
 Topology,
\textbf{38}, No. 3,  pp. 529 -- 541
(1999).

}

\newcommand{\paadv}
{  A.V.Pajitnov,
\emph{   $C^0$-generic properties of
boundary operators in the Novikov
complex },
\quad E-print:
math.DG/9812157, 29 Dec 1998,
journal article:
Advances in Mathematical Sciences,
 vol. 197, 1999, p.29 -- 117.
}

\newcommand{\pawitt}
{  A.V.Pajitnov,
\emph{Closed orbits of gradient
flows and logarithms of non-abelian Witt vectors},
\quad E-print:
 math.DG/9908010, 2 Aug. 1999
journal article:
  K-theory, Vol. 21 No. 4, 2000.
}

\newcommand{\pajandran}
{A.V.Pajitnov, A.Ranicki,
\emph{The Whitehead group of the Novikov ring},
\quad E-print:
 math.AT/0012031, 5 dec 2000,
 journal article:
  K-theory, Vol. 21 No. 4, 2000.
}

\newcommand{\paodense}
{ A.V.Pazhitnov,
\emph{
On the Novikov
complex for rational Morse forms
}, Preprint:
~
Institut for Matematik og datalogi
 Odense Universitet,
 Preprints 1991, No. 12
 Odense, October 1991.
 }

%BOOKS IN RUSSIAN

\newcommand{\postnTM}
{   M.M.Postnikov,
\emph{Introduction to Morse theory
},
 Moscow, Nauka,  1971, in Russian
}

%BOOKS IN FOREIGN LANGUAGES

\newcommand{\arnoldEquadiff}
{   V.I.Arnold,
\emph{Ordinary Differential Equations
},
 Moscow, Nauka,  1975.
}

\newcommand{\abrob}
{   R. Abraham, J.Robbin,
\emph{Transversal mappings and flows
},
 Benjamin, New York,  1967.
}

\newcommand{\atiyahandmacdo}
{   M.F.Atiyah, I.G.Macdonald
\emph{Introduction to commutative algebra
},
Addison-Wesley,   1969.
}

\newcommand{\bass}
{ H.Bass,
\emph{Algebraic K-theory},
Benjamin, 1968.
}

\newcommand{\birkrota}
{ G.Birkhoff, G-C. Rota,
\emph{Ordinary differential equations},
Blaisdell Publishing Company, 1962.
}

\newcommand{\bour}
{N.Bourbaki,
\emph{Groupes de Lie, Alg\`ebres de Lie
}
}

\newcommand{\cartanCD}
{   H.Cartan
\emph{Cours de Calcul Diff\'erentiel
},
Hermann,  1977.
}

\newcommand{\cohen}{
M.M.Cohen
\emph{A course in Simple-Homotopy theory},
Springer, 1972.}

\newcommand{\dieudonne}
{   J.Dieudonn\'e
\emph{Foundations of modern analysis
},
Academic press, 1960.
}

\newcommand{\dieud}
{   J.Dieudonn\'e
\emph{El\'ements d'analyse, Vol. III
},
Gauthier-Villars, 1970
}

\newcommand{\dold}
{   A.Dold
\emph{Lectures on Algebraic Topology
},
Springer,  1972.
}

\newcommand{\irwinB}
{ M.C.Irwin
\emph{Smooth dynamical systems
},
Academic Press,  1980.
}

\newcommand{\hirsch}
{ M.Hirsch
\emph{Differential Topology
},
Springer,  1976.
}

\newcommand{\huse}
{ D.Husemoller
\emph{Fibre bundles
},
McGraw-Hill, 1966.
}

\newcommand{\kiang}
{   Kiang Tsai-han
\emph{The theory of Fixed point classes},
 Springer, 1989.
}

\newcommand{\kling}
{   W.~Klingenberg
\emph{Lectures on closed geodesics
},
 Springer, 1978.
}

\newcommand{\liapounov}
{   Liapounov A.M.,
\emph{Probl\`eme g\'en\'eral de la stabilit\'e du mouvement},
 Princeton University
Press,
 1947.
}

\newcommand{\milnMT}
{   J.~Milnor,
\emph{ Morse theory
},
 Princeton University Press, 1963.
}

\newcommand{\milnKT}
{   J.~Milnor,
\emph{ Introduction to algebraic K-theory
},
 Princeton University Press, 1971.
}

\newcommand{\milnhcob}
{   J.~Milnor,
\emph{Lectures on the
h-cobordism theorem},
 Princeton University
Press,
 1965.
}

\newcommand{\milnstash}
{   J.Milnor and J.Stasheff,
\emph{ Characteristic Classes},
 Princeton University Press,
 1974.}

\newcommand{\morse}
{  M.Morse,
\emph{Calculus of Variations in the Large},
  American Mathematical Society Colloquium Publications,
Vol.18,
 1934.}

\newcommand{\munkres}
{  J.R.Munkres,
\emph{Elementary differential toplogy},
Annals of Math. Studies,
Vol.54, Pinceton
 1963.}

\newcommand{\cohn}
{ P.M.Cohn,
\emph{Free rings and their relations},
  Academic press
( 1971)}

\newcommand{\lam}
{    T.Y.Lam,
\emph{Serre's Conjecture,         }
Lecture Notes in Mathematics {\bff 635}, (1978) 227 p.
}

\newcommand{\lang}
{    S.Lang,
\emph{Algebra ,         }
Addison-Wesley (1965)
}

\newcommand{\massey}
{   W.Massey,
\emph{ Homology and cohomology theory
},
 Marcel Dekker, 1978.
}

\newcommand{\palisdemelo}
{   J.Palis, Jr.,
W.de Melo,
\emph{ Geometric  theory of dynamical systems},
 Springer, 1982.
}

\newcommand{\ranKL}
{   A.A.Ranicki,
\emph{Lower $K$- and $L$-theory,         }
LMS Lecture Notes 178, Cambridge, 1992
}

\newcommand{\ranKNO}
{   A.A.Ranicki,
\emph{High-dimensional knot theory,         }
Springer, 1998
}

\newcommand{\ranitor}
{   A.A.Ranicki,
\emph{The algebraic theory of torsion I.,}
Lecture Notes in Mathematics {\bff 1126} (1985), 199--237.
}

\newcommand{\rock}
{   R.T.Rockafellar,
\emph{Convex Analysis},
Princeton University Press (1970).
}

\newcommand{\rosenberg}
{   J.Rosenberg,
\emph{Algebraic $K$-theory
 and its applications},
Springer, (1994).
}

\def\sharko{ V.V.Sharko,
\emph{ Funktsii na mnogoobraziyah
(algebraicheskie i topologicheskie aspekty)}, Kiev,
Naukova dumka
(1990),
 31-35.
   \quad English translation:
 }

\newcommand{\spanier}
{ E.H.Spanier,
\emph{Algebraic topology},
  McGraw-Hill,
(1966)}

\newcommand{\stong}
{R.Stong,
\emph{Notes on Cobordism theory},
Princeton, New Jersey, 1968}

\newcommand{\switzer}
{ R.M.Switzer,
\emph{Algebraic topology -- homotopy and homology},
  Springer,
( 1975)}

\newcommand{\wander}
{ B.L.van der Waerden,
\emph{Algebra 1},
Springer, ,
( 1971)}

\newcommand{\weibel}
{ C.A.Weibel,
\emph{An introduction to homological algebra },
Cambridge University press,
( 1997)}

%%%%%%%%%%%%%%%%%%%%%%%%%%%%%%%%%%%%
%%%%%%%%%%%%%%%%%%%%%%%%%%%%%%%%%%%%
%PREPRINTS
%%%%%%%%%%%%%%%%%%%%%%%%%%%%%%%%%%%%
%%%%%%%%%%%%%%%%%%%%%%%%%%%%%%%%%%%%

%\newcommand{\huli}{   M.Hutchings, Y.J.Lee
%\emph{              Circle-valued Morse theory, Reidemeister torsion
%and Seiberg-Witten invariants of 3-manifolds},
%\quad  E-print dg-ga/9612004 3  Dec 1996.}

\newcommand{\burghelea}{ D.Burghelea, S.Haller
\emph{ On the topology and analysis of a closed 1-form
(Novikov's theory revisited)
}, \quad  E-print dg-ga/0101043  5 jan. 2001}

\newcommand{\hulihuli}{ M.Hutchings, Y.J.Lee,
\emph{ Circle-valued Morse theory and Reidemeister Torsion
}, \quad  E-print dg-ga/9706012 23 June 1997}

\newcommand{\mengt}
{ G.Meng, C.H.Taubes
\emph{SW=Milnor Torsion }
\quad   preprint
(1996)
}

\newcommand{\ranprepr}
{  A.A.Ranicki,
      \emph{
Finite domination and Novikov rings},
preprint,
 1993 }

\newcommand{\sheiham}
{  D.Sheiham,
\emph{Noncommutative characteristic
polynomials and Cohn localization},
\quad   preprint, 2000}

%%%%%%%%%%%

%%%%%%%%%%%%%%

%ii\input /home/a/bibR_E

\def\kozlovskii{O.S.Kozlovskii
\emph{The dynamics
of intersections of
analytic manifolds
}, Doklady ANSSSR,
{\bf 323}
(1992).
\quad English translation:

 Sov.Math.Dokl.
{\bff 45}
(1992), 425--427.}

\def\novidok{S.P.Novikov,
\emph{ Many-valued functions
 and functionals. An analogue of Morse theory  },
Doklady AN SSSR,
{\bf 260}
(1981),  31-35 (in Russian),
\quad English translation:
 Sov.Math.Dokl.
{\bff 24}
(1981), 222-226. }

\newcommand{\noviuspe}
{ S.P.Novikov,
\emph{The hamiltonian formalism and a
multivalued analogue of
Morse theory,
}
Uspekhi Mat. Nauk,
{\bff 37}    (1982),  3-49(in Russian),
\quad English translation:
 Russ. Math. Surveys,
{\bff 37} (1982),
 1--56.}

\newcommand{\pafest}
{  A.V.Pajitnov,
\emph{ Simple homotopy type of Novikov complex
and $\zeta$-function of the gradient flow, }
\quad E-print:
dg-ga/970614 9 July 1997;
journal article:
Russian Mathematical Surveys,
\textbf{54}
(1999), 117 -- 170.}

\newcommand{
\pastpet}
{  A.V.Pajitnov,
{\it   The incidence coefficients in the Novikov
complex are
generically rational functions,}
\quad E-print:  dg-ga/9603006 14 March  96,
journal article: Algebra i Analiz,
{\bff 9}, no.5 (1997),  92--139 (in Russian),
\quad English translation:
Sankt-Petersbourg Mathematical Journal
\textbf{9}
(1998),
no. 5, p. 969 -- 1006.
}

\newcommand{\paclo}
{  A.V.Pajitnov,
\emph{ Counting closed orbits
of gradients of circle-valued maps},~
E-print:  math.DG/0104273 28 Apr. 2001,
journal article:
Algebra i Analiz,
{\bff 14}, no.3 (2002),  92--139
(in Russian),
English translation:
Sankt-Petersbourg Mathematical Journal.
{\bff 14}, no.3 (2003).
}
%%%%%%%%%%%%%%%%%%%%%%%%%%%%%%%
%%%%%%%%%%%%%%%%%%%%%%%%%%%%%%%%
%%%%%%%%%%%%%%%%%%%%%%%%%%%%%%%%%%

\title{$C^0$-topology
in Morse theory}
\author{A.V.Pajitnov}
\address{Laboratoire Math\'ematiques Jean Leray UMR 6629,
Universit\'e de Nantes,
Facult\'e des Sciences,
2, rue de la Houssini\`ere,
44072, Nantes, Cedex}
\email{ pajitnov@math.univ-nantes.fr}
%ii\input abstract
\begin{abstract}
Let $f$ be a Morse function
on a closed manifold
$M$, and $v$ be a Riemannian   gradient
of $f$ satisfying the transversality condition.
 The classical construction
(due to Morse, Smale, Thom, Witten),
based on the counting of flow lines joining critical points
of the function $f$ associates  to these data the
Morse complex $\MM_*(f,v)$.

In the present paper
we introduce  a new class
of vector fields ({\it $f$-gradients})
associated to a Morse function $f$.
This class is wider than the class of
Riemannian gradients and
provides a natural framework for the
study of the Morse complex.
Our  construction of the Morse complex does not use
the counting of the flow lines, but rather
the fundamental classes of the stable manifolds,
and this allows to replace the transversality
condition required in the classical setting
by a weaker condition on the $f$-gradient
({\it almost transversality condition})
which is $C^0$-stable.
We  prove then that the Morse  complex is stable
\wrt~  $C^0$-small perturbations of the $f$-gradient,
and study  the functorial properties of the Morse complex.

The last two sections of the
paper are devoted to
the properties of functoriality
and $C^0$-stability
for the Novikov complex $\NN_*(f,v)$ where $f$ is
 a  circle-valued Morse map
and $v$ is an almost transverse $f$-gradient.
\end{abstract}
\maketitle

%ii\input intro
\section{Introduction}
\mlb{s:intro}

Recall that the  Riemannian gradient of  a
differentiable
function $f:M\to\RRR$ on a Riemannian
manifold $M$ is defined by the following
formula:
$$
\langle\grad f(x),h\rangle
=
f'(x)(h)
$$
(where  $\langle\cdot,\cdot\rangle$
stands for the scalar product on  $T_xM$,
and $h\in T_xM$).
The function $f$
is strictly increasing along any
non-constant integral curve of $\grad f$.
Thus the properties of $f$ and of
the flow generated by
$\grad f$
({\it gradient flow})
are closely related
to each other.

In the Morse theory the use of
the gradient flows was initiated by
R.Thom in the article \cite{thom}.
Later on, the techniques of
gradient flows in Morse theory
were  developed
in several  papers of M.Morse
and in the book
"Lectures on the $h$-cobordism theorem"
by J.Milnor \cite{milnhcob}
which provides an alternative approach to
S.Smale's proof of $h$-cobordism conjecture.
In this book J.Milnor works with   a certain
particular class of Riemannian  gradients;
here is the definition.
\footnote{
It is well-known (and easy to prove),
  that each gradient-like vector field
is a Riemannian gradient.}
\bede[\cite{milnhcob}, \S 3]\mlb{d:glvf}
Let $M$ be a manifold, $f:M\to\RRR$
be a Morse function. A vector field $v$ is
called {\it \glvf~ }  for $f$, if

1) for every $x\notin S(f)$ we have:
$f'(x)(v(x))>0$,

2) for every $p\in S(f)$ there is a  chart
$\Psi:U\to V\sbs \RRR^m$
of the manifold $M$, \sut
\begin{gather}\lb{f:chart_vv}
(f\circ\Psi^{-1})(x_1,...,x_m)
=
f(p)-\sum_{i=1}^k x_i^2
+
\sum_{i=k+1}^m x_i^2,
 \\
\Psi_*(v)(x_1,...,x_m)
=
(-x_1,...,-x_k, x_{k+1}, ... x_m)
\quad\mx{ where }\quad k=\ind p.
\end{gather}
\end{defi}
(Here $S(f)$ stands for the set of all critical points of $f$.)
This notion has
many advantages from the
point of view of the differential topology .
To construct and modify such a vector
field there is no
need for an auxiliary  object such as  Riemannian metric.
Also the local structure of the vector
field and its integral curves
nearby the  critical points of $f$ is much simpler
than for general Riemannian gradients.

On the other hand the class of gradient-like
vector fields
is "smaller" than one would like it to be; for example,
a  small
perturbation of a gradient-like vector field for $f$
is not necessarily again a \glvf~ for $f$.

In Section \mrf{s:gra} of the present  paper   we suggest
a  class of vector fields, which
is strictly larger than the class of
\rgr s, but its definition,
similarly to the definition of \glvf s,
uses only the condition 1) above and a
certain local non-degeneracy
condition at every critical point
(Definition \mrf{d:f_grad}).
We call these vector fields {\it $f$-gradients}.
The subject of  Section \mrf{s:morcom} is the
construction of the Morse complex for a
Morse function $f$ and an $f$-gradient $v$.
Instead of counting the flow lines joining the critical points
our construction is based on  the fundamental classes of the descending discs.
This allows to weaken the  condition of transversality
required in the classical definition of the Morse complex.
Namely we use the
{\it almost transversality condition}
(see Definition \mrf{d:trans_cond}).
Intuitively, the almost transversality condition
requires that the stable and unstable discs of two
critical points  be transverse only
for the case when the sum of the dimensions of these discs
is not greater than the dimension of the ambient manifold
(so that these discs  do not intersect).
The advantage of this condition is that
the set of almost transverse $f$-gradients is
open in the set of all $f$-gradients \wrt~ $C^0$-topology,
and this leads to  a natural formulation
of the $C^0$-stability property of the
Morse complex (see Theorem \mrf{t:general_stable}).
In the same section we study the functorial properties of
the Morse complex.

The subject of Section \mrf{s:nov_com}
is the Novikov complex.
For every Morse map $f:M\to S^1$
and every almost transverse $f$-gradient $v$ we construct
a chain complex $\NN_*(f,v)$
of modules over the ring $\ZZZ((t))$
of Laurent series.
Our main aim here is to construct
a canonical chain equivalence between
this version of the Novikov complex to
the Novikov completion of the
singular chain complex of the corresponding
infinite cyclic covering.
We study the functorial properties
of this chain equivalence.

In the section \mrf{s:further}
we discuss the $C^0$-stability of the
Novikov complex, and a formula which
relates the homotopy class  of the
canonical chain equivalence
constructed in the previous section
and the Lefschetz zeta function of the gradient flow.

%ii\input gra
\section{Gradients of Morse functions and forms }
\mlb{s:gra}

\subsection{Definition of $f$-gradients}
\mlb{su:def_grad}

Let $W$ be a compact cobordism
or a manifold without boundary,
and let $f:W\to \RRR$ be a Morse function.
Let $v$ be a $\smo$vector field on $W$
satisfying the following condition
\bq
\lb{f:weak_gra}
f'(x)(v(x)) >0
\mx{~ whenever ~ }
x\notin S(f).
\end{equation}
The function
$\phi(x)=f'(x)(v(x))$
vanishes on
$S(f)$ and is strictly
positive on $W\sm S(f)$.
Therefore every point  $p\in S(f)$
is a point of local minimum of $\phi$, and
$\phi'(p)=0$.
\bede\mlb{d:f_grad}
A $\smo$ vector field $v$ is called {\it $f$-gradient}
if the condition \rrf{f:weak_gra} holds, and
every $p\in S(f)$
is a point of non-degenerate minimum of the function
$\phi(x)=f'(x)(v(x))$ (that is, the second derivative
$\phi''(p)$ is a non-degenerate bilinear form on $T_pW$).
\end{defi}
We shall now compute the second derivative $\phi''(p)$
in terms of the derivatives of $f$ and $v$
and give an alternative formulation
of the non-degeneracy condition above.

\bele\mlb{l:van_critpt}
Let $v$ be an $f$-gradient,
and $p\in S(f)$.
Then $v(p)=0$.
\enle
\Prf
Pick a Morse chart $\Psi:U\to V\sbs \RRR^n$ ~
for $f$,
that is a chart where the function $f\circ\Psi^{-1}-f(p)$
equals  the quadratic form
$$
Q(x,y)=-||x||^2+||y||^2;
\quad
x\in\RRR^k, y\in \RRR^{n-k}
\mx{ ~(where~ } k=\ind p),
$$
and let $w=\Psi_*v$.
We have $Q'(z)(w(z))>0$ for every $z\not=(0,0)$,
and we must prove $w(0,0)=0$.
Write $w(0,0)=(\xi, \eta)$, and assume that
$\xi\not=0$.
Write
$$ w(x,0)=(\xi+\OO_1(x), \eta+\OO_2(x))$$
where
$||\OO_i(x)||\leq C||x||$
nearby $0$.
We have
\bq\lb{f:posit}
Q'(x,0)(w(x,0))
=
-\langle 2x, \xi +\OO_1(x) \rangle
\geq 0
\mxx{ for every } x.
\end{equation}
Set $x=t\xi$ with $t\to 0$
to deduce  that \rrf{f:posit} is possible only for $\xi=0$.
Similarly we prove that $\eta=0$.
$\qs$

At every point $p\in S(f)$
the vector field $v$ vanishes, and therefore
there is  a well defined linear map
 $v'(p):T_pW\to T_pW$.
A simple computation in the local coordinates proves
the following lemma.
\bele\mlb{l:form_sec_der}
For $p\in S(f)$ we have
\bq\lb{f:form_sec_der}
\phi''(p)(h,k)
=
f''(p)(v'(p)h,k)
+
f''(p)(v'(p)k,h).\hspace{2cm}\qs
\end{equation}

\enle

The next proposition follows immediately.
\bepr\mlb{p:f_grad}
A vector field $v$ on $W$
is an  {\it $f$-gradient}
\ifff~

A) for every $x\notin S(f)$
we have
$f'(x)(v(x))>0$ and

B)  for every $ p\in S(f)$ we have
\bq\lb{f:def_grad}
f''(p)(v'(p)h,h)>0 \mx{ for every }
h\in T_pW, ~
h\not=0.\hspace{2cm}\qs
\end{equation}
\enpr

Both \glvf s and Riemannian gradients
are examples of $f$-gradients. This
is quite obvious for
the \glvf s and for the Riemannian gradients it is the
subject of the next proposition.

\bepr\mlb{p:riem_ok}
Every Riemannian gradient
is an $f$-gradient.
\enpr
\Prf The
condition A) of the definition
 is obviously satisfied.
To prove  the condition B)
it suffices to consider the case when
$f$ is defined in a \nei~~ of $0$ in $\RRR^n$,
and $0$ is the only critical point of $f$ in this \nei.
Let $\GG(x)$ be   the matrix
of the Riemannian metric \wrt~
the coordinates in the $\RRR^m$.
Let $\nabla f(x)$ denote the Euclidean gradient
$(\chape {f}{x_1}(x),..., \chape {f}{x_m}(x))$
of $f$.
We can assume that in a \nei~ of $0$ we have:
$$
\GG(x) =\Id+
\OO_1(||x||^2)
$$
where  $\OO_1(||x||^2)\leq C||x||^2$.
Applying a linear change of base in
$\RRR^n$, we can assume also that
$$
\nabla f(x)= Ax+\OO_2(||x||^2),
$$
where $A$ is a diagonal matrix with
non-zero diagonal entries.
The Riemannian gradient of $f$
\wrt~ our metric
satisfies
$$
\grad f(x) = \GG^{-1}(x)\cdot \nabla f(x) = Ax+\OO_3(||x||^2).
$$
Therefore
\bq\lb{f:grad_formu}
f'(x)(\grad f(x))=\langle \grad f(x),\grad f(x) \rangle
=
\langle Ax, Ax \rangle +\OO(||x||^3)
\end{equation}
and this function has a non-degenerate
minimum at $x=0$.$\qs$

\subsection{Topological properties
of the space of all $f$-gradients}
\mlb{su:space_grad}

In this subsection $W$ is a cobordism, and
$f:W\to[a,b]$ is a Morse function.
The space of all $\smo$ vector fields on $W$
endowed with the usual $\smo$  topology
will be denoted $\VV(W)$.
In this space consider the subspace
$\VV(f)$
of all the vector fields, which vanish on $S(f)$.
This is a  closed subspace of finite codimension.
The space $G(f)$ of all $f$-gradients
 is a subspace of $\VV(f)$.

\bepr\mlb{p:open_grad}
The set $G(f)$ is an
open convex subset of
$\VV(f)$.
\enpr
\Prf
Convexity is  obvious. As for  the
openness, observe that the  condition
B) from the definition
 of $f$-gradients is clearly open
\wrt~$\smo$ topology.
Proceed to the condition A).
Let $w\in G(f)$. We have to prove that a \nei~
of $w$ in $\VV(f)$ is in $G(f)$.
The cobordism $W$ is compact, therefore
it suffices to prove that
for every $a\in W$
the following property holds:
\begin{quote}
$(\PP)$\quad
There is a \nei~ $U(a)$ of $a$
and  a \nei~ $\VV$ of $w$ in $\VV(f)$
such that for every $u\in\VV$ and every
$x\in U(a)\sm S(f)$
we have:
\bq\lb{f:pluss}
f'(x)(u(x))>0.
\end{equation}
\end{quote}

This property
is obviously fulfilled  for any
point $a\notin S(f)$.
In  the case $a\in S(f)$
the property $(\PP)$
follows from the next
Lemma,
which is proved by an easy application of the
Taylor development.
For a function $g:U\to E$
defined in an open set $U\sbs \RRR^n$,
with values in a normed space $E$, put
$$
||g||_U
=
\sup_{x\in U} ||g(x)||.
$$
\bele\mlb{l:open_loc}
Let $\phi:U\to\RRR$
be a $\smo$ function defined in a
\nei~ $U$ of $0$ in $\RRR^m$.
Assume that $\phi(0)=0$ and that
$\phi$ has a non-degenerate minimum in $0$.
Then there is a compact \nei~ $U$ of $0$ and
a number $\d>0$ \sut~
every $\smo$ function $\psi:U\to\RRR$
satisfying
\bq\lb{f:delta}
\psi(0)=0,~ \psi'(0)=0, ~
||\phi''-\psi''||_U<\d,~
||\phi'''-\psi'''||_U<\d,~
\end{equation}
satisfies also $\psi(x)>0$ if $\in U $ and
$x\not=0$. $\qs$
\enle

Now we have already three types of
vector fields, associated with a given Morse function:
Riemannian gradients, gradient-like
vector fields, and $f$-gradients.
One more notion is useful:

\bede\mlb{d:weak_gra}
A vector field $v$ is called {\it a weak gradient for $f$}
\footnote{The notion of weak gradient was
introduced by D.Sch\"utz
in his paper \cite{schuetzCI}.}
if
$$
f'(x)(v(x))>0
\mxx{ for every }
x\notin S(f).
$$
\end{defi}
The set of all weak
gradients for $f$ will be denoted by
$GW(f)$.
We have  the  inclusions
$$
GL(f)\sbs GR(f)\sbs G(f)\sbs GW(f)\sbs \VV(f),
$$
and the rest of this section is devoted to the study
of their topological properties.
\bepr\mlb{p:grad_dense}
The set $G(f)$ is everywhere
dense in $GW(f)$.
\enpr
\Prf
Let $v$ be any weak
gradient for $f$, and
$v_0$ be any $f$-gradient
(for example a \glvf~ for $f$).
The vector field
 $w=v+\e v_0$
is an $f$-gradient, which converges to $v$ as $\e\to 0$.
$\qs$

Now let us compare  the spaces
$GR(f), GL(f)$ and $G(f)$.
The reader has certainly noticed that
the main difference between
$f$-gradients,  Riemannian gradients
and \glvf s
is in their behavior nearby the critical points.
It is easy to show that
for a critical point $p\in S(f)$
the linear map $v'(p):T_pW\to T_pW$
is
1) diagonalizable over $\RRR$
if $v$ is a Riemannian gradient
2) has only the eigenvalues $\pm 1$
if $v$ is a \glvf.
Using this observation it is not difficult to show that
for every Morse function $f$ with $S(f)\not= 0$,
the space $GL(f)$ is not dense in $GR(f)$, and if moreover,
$\dim W >2$, the space $GR(f)$ is not dense in
$G(f)$.

\subsection{Gradients of Morse forms}
\mlb{su:gra_morse_forms}

The techniques of the previous
subsections are purely local, even when the
statements of the results
are of global character. Thus the contents
of these subsections is generalized  to
the case of Morse forms without any difficulty.
In this subsection we list  the corresponding
definitions and results.
We shall consider  manifolds
without boundary, since only this case will be used in
the later applications. The most important for us is the
case when the Morse form $\o$ is of the form
$\o=df$, where $f:M\to S^1$ is a Morse map.
Let $M$ be a \ma~ without boundary, $\o$
be a Morse form on $M$.

\bede\mlb{d:o_grad}
We say that  $v$
is an  {\it $\o$-gradient}
if $\o(x)(v(x))>0$
whenever $x\notin S(f)$
 and  for every $ p\in S(f)$
the real-valued function $\phi(x)=\o(x)(v(x))$
has a non-degenerate local minimum at $p$.
The space of all $\o$-gradients
will be denoted $G(\o)$.
\end{defi}
Similarly to the case of the Morse
 functions, the condition (B)
is equivalent to the following condition:
\bq\lb{f:nondeg_forms}
f''(p)(v'(p)h,h)>0 \mxx{ for every } h\not=0, h\in T_pM
\end{equation}
(where $f$ is any function in a \nei~ of $p$ with $df=\o$).

\bede\mlb{d:riem_grad_forms}
The {\it \rgr~ } $\grad \o$
\wrt~ a given Riemannian metric~ on $M$
is defined by
\bq\lb{f:riemgr_forms}
\o(x)(h)
=
\langle
\grad\o(x),
h
\rangle
\end{equation}
where $x\in M, ~ h\in T_xM$.
The space of all  Riemannian gradients for $\o$
will be denoted $GR(\o)$.
\end{defi}
A Riemannian gradient for $\o$ is an $\o$-gradient
(similarly to Proposition  \mrf{p:riem_ok}).

\bede\mlb{d:glvf_forms}
A vector field $v$
is called {\it \glvf~ for $\o$},
if for open set $U\sbs M$
and every Morse function $f:M\to\RRR$
with $df=\o|U$
the vector field $v|U$ is a \glvf~ for $f|U$.
The space of all \glvf s for $\o$
will be denoted $GL(\o)$.
\end{defi}
\bede\mlb{d:weak_forms}
A vector field $v$ on $M$ is called
{\it weak gradient for $\o$} if
\bq
\lb{f:weak_gra_forms}
\o(x)(v(x)) >0
\mx{~ whenever ~ }
x\notin S(\o)
\end{equation}
The space of all \glvf s for $\o$
will be denoted $GW(\o)$.
\end{defi}

Denote by $\VV(\o)$
the space of all vector fields on $M$ vanishing
on $S(\o)$. This is a closed vector subspace of
codimension $m(\o)$ in $\VV(M)$.
We have the inclusions
$$
GL(\o)\sbs GR(\o)\sbs G(\o)\sbs GW(\o)\sbs \VV(\o),
$$

The proofs of the next propostions
are completely similar to the proofs
of Propositions \mrf{p:grad_dense}
and \mrf{p:open_grad}.

\bepr\mlb{p:open_dense}
Let $W$ be a closed \ma, and
$\o$ be a closed 1-form on $W$.
The set $G(\o)$  is a
convex open subset of $\VV(\o)$.
The set $G(\o)$ is dense in the space
$GW(\o)$. $\qs$
\enpr

\bepr\mlb{p:not_dense}
If $\o$ has at least one zero
and $\dim M >2$,  then the
set $GL(\o)$ is not dense in $GR(\o)$,
and the set $GR(\o)$ is not dense in $G(\o)$.$\qs$
\enpr

\subsection{Transversality properties}
\mlb{su:transv}

Let $v$ be a vector field on a manifold $M$,
and let $v(p)=0$. Recall that the {\it stable set
$W^{st}(p,v)$ of $p$} is the set of all points
$x\in M$, such that the $v$-trajectory $\g(x, \cdot;v)$
is defined for all $t\geq 0$ and converges to
$p$ as $t\to\infty$.
The {\it unstable set} $W^{un}(p,v)$ is by definition
the set $W^{st}(p,-v)$.

Let $M$ be a \ma~ without boundary, $\o$ be a closed
1-form, $v$ be an $f$-gradient,
$p$ be a zero of $f$.
It is easy to deduce form the formula \rrf{f:nondeg_forms}
that $p$ is an elementary zero of $v$, that
is, the linear map $v'(p)$ has no imaginary eigenvalues,
and therefore by Hadamard-Perron theorem
(see \cite{abrob}, Th. 27)
there exist local stable and unstable manifolds for $v$
at $p$. In  the case when   $M$ is compact,
the (global ) stable manifold $W^{st}(p,v)$
is an immersed manifold of dimension $\ind p$,
and $W^{un}(p,v)$
is an immersed manifold of dimension $n-\ind p$.
The stable manifold $W^{st}(p,v)$
will be also called {\it descending disc}
of $p$ and denoted $D(p,v)$.
The stable manifold $W^{st}(p,-v)$
will be also called {\it ascending disc}
of $p$ and denoted $D(p,-v)$.

In the case when $v$ is an $f$-gradient,
where $f:W\to[a,b]$
is a Morse function on a cobordism,
the stable manifold $W^{st}(p,v)$
of a critical point $p\in S(f)$
is a submanifold with boundary of $W$ of
dimension $\ind p$, and $W^{st}(p,v)$
is an $(n-\ind p)$-dimensional
\sma~with boundary.

In the next definition $v$ is an $f$-gradient where
$f:W\to[a,b]$
is a Morse function, or an $\o$-gradient,
where $\o$ is a Morse form on a closed \ma.

\begin{defi}
\mlb{d:trans_cond} We say that
 $v$ is {\it transverse } or
           {\it satisfies transversality condition,} if
\bq\lb{f:trans_cond}
\big(
x,y \in S(f) \big)
\Rightarrow \big( D(x,v) \pitchfork D(y,-v) \big).
\end{equation}

\vskip0.1in

We say that
$v$ is {\it almost
transverse} or
 {\it satisfies Almost Transversality Condition},
if
\bq\lb{f:alm_trans_cond}
\big( x,y \in S(f) ~\&~
\ind x \leq \ind y \big) \Rightarrow \big(
D(x,v) \pitchfork D(y,-v)\big).
\end{equation}
\vskip0.1in
\end{defi}
The notion of almost transverse gradients
was introduced in
\cite{pasur}.
We have the corresponding version
of Kupka-Smale theorem
(the proof is similar to the classical proof,
see \cite{peixoto}, \cite{palisdemelo}).
\beth[Kupka-Smale theorem]
Let $\o$ be a closed 1-form on a
closed manifold $M$.
The set of all transverse $f$-gradients
is residual in $\VV(\o)$
(\wrt~ $\smo$ topology).
\enth

The case of Morse functions on
cobordisms is not covered
by this theorem, but it can be
dealt with by a similar  argument.
\bepr
Let $f:W\to[a,b]$ be a Morse function.
The set of all transverse $f$-gradients
is residual in $\VV(f)$
(\wrt~ $\smo$ topology).
\enpr

The set of all almost transverse gradients
is much larger than the set of transverse
gradients, as the next proposition shows
(the proof is similar to the proof of
Cor. 1.7 of \cite{pastpet}).
\bepr
Let $f:W\to[a,b]$ be a Morse function
on  a cobordism $W$.
The set of all almost
transverse $f$-gradients
is dense  in $\VV(f)$
\wrt~ $\smo$ topology,
and open in $\VV(f)$ \wrt~ $C^0$-topology.
\enpr

\subsection{Ordered Morse functions and
Morse-Smale filtrations}
\mlb{su:morse_smale}

This section contains a bunch of definitions
to be used in the sequel.

\begin{defi}
{\it An ordered Morse function}
is a Morse function $f:W\to[a,b]$
together with a sequence
$a_0,...,a_{n+1}$
of regular values such that~
$a=a_0<a_1<...<a_{n+1}=b$
and
 \begin{equation}
S_i(f)\subset
 f^{-1}(]a_i,a_{i+1}[) \lb{f:ordered}
\end{equation}
for every $i$.
The  sequence
$a_0,..., a_{n+1}$
satisfying  (\ref{f:ordered})
is called {\it ordering sequence} for $f$.
\end{defi}

Note that we consider the ordering
sequence as a part
of the data of the ordered Morse function.
Each ordered Morse function generates
a filtration of $W$, defined by
\bq\lb{f:ms_filtrr}
W^{(k)}=\phi^{-1}([a_0, a_{k+1}]);
\end{equation}

\bede\mlb{d:morse_smale_filtr}
A filtration
$W^{(k)}$
of the cobordism $W$
(where $0\leq k\leq \dim W$) is
called {\it Morse-Smale filtration}
if there is an ordered  Morse function
 $\phi:W\to [a,b]$, such that
\rrf{f:ms_filtrr}
holds.
\end{defi}

\begin{defi}
\mlb{d:adjust}
A Morse function $\phi:W\to[\a,\beta]$
is
called {\it adjusted to $(f,v)$}, if:

 1) $S(\phi)=S(f)$,

 2) the function $f-\phi$ is constant in a \nei~
of $\pr_0 W$, in a \nei~ of $\pr_1 W$, and in a \nei~ of each
point of
$S(f)$,

3) $v$ is also a $\phi$-\gr.
\end{defi}

An  analog of the classical Rearrangement Lemma
(see \cite{milnhcob}, \S 4, or \cite{paadv}, Prop. 2.23)
is also valid in the context of $f$-gradients, and this
implies the next proposition.
\bepr\mlb{p:adj}
Let $\fcob$ be  a Morse function on a
cobordism $W$, and $v$ be an
almost transverse $f$-gradient. There
is an ordered Morse function $\phi:W\to [a,b]$,
adjusted to $(f,v)$.
such that $v$ is a $\phi$-gradient.$\qs$
\enpr

The next proposition is proved by the same argument as
Lemma 2.74 of \cite{paadv}.
Here and elsewhere in this paper the symbol $||\cdot ||$
denotes the {\it $C^0$-norm}.
\bepr\mlb{p:adj_perturb}
Let $g$ be a Morse function on $W$, adjusted to $(f,v)$.
Then there is $\d>0$ such that for every $f$-gradient
$w$ with $||w-v||<\d$
 the function $g$ is adjusted to $(f,w)$. $\qs$
\enpr

\bede\mlb{d:morse_sm_adj}
A Morse-Smale filtration
$W^{(k)}$ on $W$
is called {\it adjusted to $(f,v)$}
if there is an ordered Morse function
$\phi$, adjusted to $(f,v)$,
such that
\rrf{f:ms_filtrr} holds.
\end{defi}

\beco\mlb{c:m_s_adj_perturb}
Let $W^{(k)}$ be a Morse-Smale
filtration of $W$, adjusted to $(f,v)$.
Then there is $\d>0$ such that for every $f$-gradient
$w$ with $||w-v||<\d$
the filtration  $W^{(k)}$
is adjusted to $(f,w)$.
\enco

%ii\input morcom

\section{Morse complex}
\mlb{s:morcom}

\mfcobv
If $v$ is  transverse, then
the classical construction of the Morse complex,
 based on  counting
of the flow lines joining critical points
(see \cite{witten}), carries over to
the present case without any changes.
Our aim in this section is to generalize this construction
to the case of {\it almost transverse }
$f$-gradients, and to
study its properties.
An $f$-gradient $v$ is called {\it oriented}
if for every critical point $p\in S(f)$
an orientation of the descending disc $D(p,v)$
is fixed.
In the first subsection we construct the Morse complex for
oriented almost transverse gradients of
 ordered Morse functions.
In the second subsection
we generalize this construction to the case
of arbitrary Morse functions and \orial~
gradients. We show that the resulting complex is
stable \wrt~ $C^0$-small perturbations of the gradient.
 We construct a canonical
chain equivalence between the Morse complex and
the singular chain complex.
In the last subsection we study the functorial
properties of these constructions.

\subsection{Ordered Morse functions}
\mlb{su:ordered_morcom}

Let $W$ be a cobordism
of dimension $n$.
Let $\phi:W\to [a,b]$
be an ordered Morse function,
 $a=a_0,..., a_{n+1}=b$ be the ordering  sequence
of $\phi$, and
$W^{(k)}=\phi^{-1}([a_0, a_{k+1}])$
be the corresponding filtration of $W$.
Let $v$ be an almost
transverse $\phi$-gradient.
As in the standard Morse theory,
this filtration is {\it cellular}, that is,
the homology $H_*(W^{(k+1)},W^{(k)})$
vanishes for every $*\not=k$, and
 $C_k=H_k(W^{(k+1)},W^{(k)})$
is a free abelian group generated by
the homology classes
of the descending discs of the critical
points of index $k$.
The generator, corresponding to the
descending disc $D(p,v)$ of a critical point $p\in S_k(f)$
with the chosen orientation
 will be denoted $\D(p,v)$.
Let us  endow the graded abelian group
$C_*$ with the boundary
operator induced from the exact sequences of  triples
$(W^{(k)},W^{(k-1)}, W^{(k-2)})$; denote the
resulting chain complex by
$C_*(\phi,v)$. This is a chain complex of
free abelian groups
endowed with a basis $\{\D(p,v)\}$.
The next theorem is an immediate consequence
of the general properties of cellular filtrations
(see \cite{dold} Ch. 5).
\bepr\mlb{p:ordered_morcom}
$H_*(C_*(\phi,v))\approx H_*(W, \dow)$.$\qs$
\enpr

The aim of the next proposition is
to compare the chain complexes associated with
an ordered Morse function $\phi$
and different $\phi$-gradients.
\bede\mlb{d:sim_orient}
Let $v,w$ be oriented $\phi$-gradients.
Their orientations are called {\it similar}
if for every $p\in S(f)$
the fundamental classes of
the manifolds
$D(p,v)$ and $D(p,w)$
induce the same generator of the group
$H_k(W_p, W_p\sm \{p\})$
(where $k=\ind p$, and $W_p=\phi^{-1}([a,F(p)])$~).
\end{defi}

Recall that $||\cdot||$ stands for the $C^0$-norm.
\bepr\mlb{p:ordered_stable}
There is $\d>0$ such that for every $\phi$-gradient $w$
with $||w-v||<\d$ the chain complexes
$C_*(\phi,v)$ and $C_*(\phi,w)$
are basis-preserving isomorphic
if the orientations of $v$ and $w$ are similar.\enpr
\Prf
If all the critical points of $f$ of a given index
have the same critical value, then
it is not difficult to show that
$\D_p(v)=\D_p(w)$ for every oriented $f$-gradient $w$
endowed with  orientation similar to the orientation of $v$.
The general case is easily reduced to this particular one
by the Rearrangement Lemma
and Proposition \mrf{p:adj_perturb}. $\qs$

Now we shall construct  a canonical
chain homotopy equivalence
between $C_*(\phi,v)$  and the singular chain complex
$\SS_*(W, \dow)$. Let $\SS_*^{(k)}=\SS_*(W^{(k+1)},W^{(k)})$
denote the filtration
induced in $\SS_*(W,\dow)$ by the subsets $W^{(k)}$.
Let us say that a chain map
$\EE:C_*(\phi,v)\to \SS_*(W,\dow)$
{\it conserves filtrations}, if
$\EE(C_k(\phi,v))\sbs\SS_k^{(k)}$ for every $k$.
A chain map $\EE$ conserving filtrations
induces for every $k$ a homomorphism
$$\EE^{adj}_{(k)}:C_k(\phi,v)\to
H_k(\SS_*^{(k)}, \SS_*^{(k-1)}).$$

\beth\mlb{t:morse_equiv_alm}
There is a
chain homotopy equivalence
$
\EE_*:C_*(\phi,v)\to \SS_*$
which
conserves filtrations and satisfies
$\EE^{adj}_{(*)}=\id $.
Such homotopy equivalence is homotopy unique.
\enth
We shall omit the proof, which repeats almost verbatim
the proof of Theorem A.5' of \cite{patou}. $\qs$
\bede
This chain equivalence will be denoted $\chi(\phi,v)$.
\end{defi}

\subsection{Morse complexes and $C^0$-stability}
\mlb{su:gen_morcom}

Now we proceed to Morse functions which are not
ordered in general.
Let $\fcob$ be a Morse function on a cobordism $W$
and $v$ be an almost tranverse oriented $f$-gradient.
Pick an ordered Morse function $\phi:W\to [a,b]$,
such that $v$ is also a $\phi$-gradient.
As the next proposition shows,
the chain complex $C_*(\phi,v)$
associated to $\phi$ and $v$
does not depend on the particular
choice of $\phi$,
neither does the homotopy class of
the corresponding
chain equivalence $\chi(\phi,v)$.

\beth\mlb{t:independ_ord}
A. Let $\phi_1, \phi_2$ be
two ordered Morse functions,
such that
$v$ is a gradient for both.
There is a homotopy commutative diagram
$$
\xymatrix{
C_*(\phi_1,v) \ar[rd]_{\chi(\phi_1,v)}
\ar[rr]^{I(\phi_1,\phi_2)} && C_*(\phi_2,v)\ar[dl]^{\chi(\phi_2,v)} \\
& \SS_*(W, \dow)
}$$
where  $I(\phi_1,\phi_2)$
is a basis preserving isomorphism of chain complexes.
\enth
\Prf
Let $W^{(k)}_1, W^{(k)}_2$
be the filtrations corresponding to the ordered
Morse functions $\phi_1, \phi_2$. In the case when
$W^{(k)}_1\sbs  W^{(k)}_2$ our assertion is easily obtained
by  functoriality.
The general case is obtained from this  particular one, since
for every two ordered Morse function $\phi_1, \phi_2$
such that $v$ is
a gradient for both, there is an  ordered function $\phi_3$,
such that $W^{(k)}_1\supset  W^{(k)}_3\sbs W^{(k)}_2$
for every $k$ and $v$ is  a $\phi_3$-gradient
(\cite{paclo}, Proposition 3.3). $\qs$

\bede\mlb{d:inc_coeff}
Let $\fcob$ be any Morse function, and $v$ be
an almost transverse $f$-gradient.
Pick any ordered Morse function
$\phi$, such that $v$ is a $\phi$-gradient.
Write
$$
\pr \big(\D(p,v) \big)=
\sum_{q\in S_{k-1}(\phi)} n(p,q;v)\cdot\D(q,v)
$$
where $\pr$ is the boundary operator in the chain complex
$C_*(\phi,v)$.
It follows from Theorem
\mrf{t:independ_ord}
that the integer  $n(p,q;v)$
does not depend on the particular choice of
an ordered function $\phi$.
This number will be called {\it incidence coefficient}
corresponding to $p,q$ and $v$.
\end{defi}

Now we can give  the definition of Morse complex
for any Morse function $\fcob$ and
its gradient $v$ satisfying
\ata.
\bede\mlb{d:inc_coef_alm}
Let $\fcob$ and $v$ be an almost transverse
oriented $f$-gradient.
Let $\MM_k(f,v)$ be the free abelian
group generated by the critical points of
$f$ of index $k$.
Define a homomorphism
$\pr_k:\MM_k(f,v)\to\MM_{k-1}(f,v)$
setting
$$
\pr_k(p)=\sum_{q\in S_{k-1}(f)} n(p,q;v) \cdot q,
$$
then $\pr_k\circ\pr_{k-1}=0$.
The resulting chain complex $\MM_*(f,v)$ is called
{\it Morse complex } corresponding to $(f,v)$.
\end{defi}
Thus the Morse complex
$\MM_*(f,v)$
is basis-preserving isomorphic to
$C_*(\phi,v)$, where $\phi$ is any
ordered Morse function such that
$v$ is a $\phi$-gradient.
It is clear that in the case when $v$ is transverse, the
Morse complex $\MM_*(f,v)$ equals the Morse
complex defined via counting
the flow lines joining the
critical points.

The next theorem follows immediately from
\mrf{p:ordered_stable}.
\beth\mlb{t:general_stable}
\mfcob.
Let $v$ be an \orial~
$f$-gradient.
Then there is $\d>0$ such that
for every $f$-gradient $w$
with $||w-v||<\d$ the chain complexes
$\MM_*(f,v)$ and $\MM_*(f,w)$
are basis-preserving isomorphic
if the orientations of $v$ and $w$ are similar.$\qs$
\enth

For every ordered Morse function $\phi:W\to\RRR$,
such that $v$ is also a $\phi$-gradient
we have the canonical chain homotopy equivalence
$\chi(\phi,v):C_*(\phi,v)\to \SS_*(W, \dow)$.
Composing it with the basis preserving isomorphism
$\MM_*(f,v)\to C_*(\phi,v)$
we obtain a  chain equivalence
\bq\lb{f:canon_equiv}
\EE=\EE(f,v):\MM_*(f,v)\rTo^\sim \SS_*(W, \dow);
\end{equation}
its homotopy class does not depend on the
particular choice of an ordered Morse function.
It  will be called {\it the Morse chain equivalence}.

\subsection{Functorial properties}
\mlb{su:functo_morse}

The aim of this subsection is to study
functorial properties of the Morse complex
and the Morse chain equivalence.
Let $f_1:M_1\to\RRR, f_2:M_2\to\RRR$
be Morse functions on closed \ma s,
and $v_1, v_2$ be \orial~ gradients for $f_1$, resp. $f_2$.
Let $A:M_1\to M_2$
be a continuous map, which satisfies the following
condition:
\begin{multline}\lb{f:alm_tr_prop}
A(D(p,v_1))~\cap~ D(q,-v_2)=\ems \\
\mxx{for every } p\in S(f_1), q\in S(f_2)
\mxx{ with} \ind p< \ind q.
\end{multline}
It is clear that
the set $\CC(v_1,v_2)$
of the  maps
satisfying the condition \rrf{f:alm_tr_prop}
is open in the space $C^0(M_1,M_2)$
of all the continuous maps
(endowed with $C^0$-topology).
As we shall see later on in this subsection the
set $\CC(v_1, v_2)$
is also dense in $C^0(M_1,M_2)$.
Assuming that the condition \rrf{f:alm_tr_prop}
is fulfilled, we shall now  construct
a homomorphism
\bq
A_\sharp:\MM_*(f_1, v_1)\to \MM_*(f_2, v_2).
\end{equation}
Pick some Morse-Smale filtrations $M_1^{(k)}, M_2^{(k)}$
of $M_1, M_2$, adjusted to $(f_1,v_1)$,
resp. to $(f_2,v_2)$.
Let us first assume that for some $T\geq 0$ the
following condition holds:
\bq\lb{f:special_ordered}
(\Phi(T,-v_2)\circ A) \big( M_1^{(k)}\big)\sbs  M_2^{(k)}
\mxx{ for every } k.
\end{equation}
(Here $\Phi(\cdot ,-v_2)$ stands for the
flow generated by $-v_2$.)
In this case the continuous map
$\Phi(T,-v_2)\circ A$
preserves the Morse-Smale filtrations, and
induces a chain map
\bq
C_*(\phi_1, v_1)
\to
C_*(\phi_2, v_2)
\end{equation}
in the Morse complexes.
It is obvious that this chain map does not depend
on the particular choice of
$T$.
We obtain therefore a chain  map
\bq\lb{f:a_sharp_phi}
A_\sharp(\phi_1, \phi_2):
\MM_*(f_1, v_1)\to \MM_*(f_2, v_2).
\end{equation}
An argument, similar to the proof of
Theorem \mrf{t:independ_ord}
shows that the map \rrf{f:a_sharp_phi}
is  independent of the choice
of $\phi_1, \phi_2, $
satisfying the condition
\rrf{f:special_ordered}.

Observe that for every continuous map $A$ satisfying
\rrf{f:alm_tr_prop}
there are Morse-Smale filtrations adjusted to $(f_i,v_i)$
such that
\rrf{f:special_ordered}
holds. Indeed, let us begin with any Morse-Smale filtration
$M_2^{(k)}$, adjusted to $(f_2, v_2)$.
Let $D_k(v_1)$ denote the union of all descending
discs  $D(p,v_1)$ with $\ind p\leq k$. The condition
\rrf{f:alm_tr_prop}
implies that for every $T$ sufficiently large,
we have
\bq\lb{f:special_orderedd}
(\Phi(T,-v_2)\circ A) \big(D_k(v_1)\big)\sbs  M_2^{(k)}
\mxx{ for every } k.
\end{equation}
Then it suffices to choose the Morse-Smale filtration of $M_1$
with
$$
M_1^{(k)}
\sbs
(\Phi(T,-v_2)\circ A)^{-1}( M_2^{(k)}) \mxx{ for every } k,
$$
(which is possible by
\cite{paclo}, Prop 3.3)
and the condition
\rrf{f:special_ordered}
is verified.

Thus for every continuous map $A$, satisfying
\rrf{f:alm_tr_prop}
we have constructed a chain map
$$
A_\sharp:
\MM_*(f_1, v_1)\to \MM_*(f_2, v_2).
$$
Our next aim is to show that this map
commutes with the Morse chain equivalences.
\bepr\mlb{p:induced_and_equiv}
The following diagram is chain homotopy commutative:
\bq\lb{f:quad}
\xymatrix{
\MM_*(f_1, v_1)\ar[d]^{\EE_1}\ar[r]^{A_\sharp} &
\MM_*(f_2, v_2)\ar[d]^{\EE_2}\\
\SS_*(M_1)\ar[r]^{A_*} & \SS_*(M_2)
}
\end{equation}
(where $A_*$ is the chain map
in the singular chain complexes induced
by $A$, and $\EE_i$ are the Morse chain equivalences).
\enpr
\Prf
Pick
any  ordered Morse functions $\phi_i$
and a positive number $T$ such that
\rrf{f:special_ordered}
holds.
In the following diagram
we denote by $\a$
the adjoint map induced by the continuous map
$A'=(\Phi(T,-v_2)\circ A)$
preserving filtrations,
and by $A'_*$ the chain map induced by $A'$
in the singular chain complexes.
This diagram
is homotopy commutative,
by  functoriality  of
the canonical chain equivalences (see \cite{patou}, Corollary 3.4
or \cite{paclo}, Corollary 2.6):
\bq\lb{f:quadd}
\xymatrix{
C_*(\phi_1, v_1)\ar[d]^{\chi(\phi_1, v_1)}\ar[r]^\a
&
C_*(\phi_2, v_2)\ar[d]^{\chi(\phi_2, v_2)}\\
\SS_*(M_1)\ar[r]^{A'_*} & \SS_*(M_2)
}
\end{equation}
It remains to note that the chain maps $A'_*$ and $A_*$
are chain homotopic. $\qs$

In the case when $A$ is a $\smo$ map,
satisfying the   condition \rrf{f:tr_prop} below,
we can obtain an explicit formula for $A_\sharp$ in terms
of intersection indices.
\begin{multline}\lb{f:tr_prop}
A|D(p,v_1):D(p,v_1)\to M_2
\mx{
is transverse to } D(q,-v_2)\\
\mx{
for every } p\in S(f_1), q\in S(f_2).
\end{multline}
A routine  argument from transversality theory
shows that the set of
all maps satisfying the above condition, is residual
in $C^\infty(M_1, M_2)$
(which implies in particular, that this set is dense,
since $C^\infty(M_1, M_2)$ is a Baire space).
Observe that the condition \rrf{f:tr_prop}
implies \rrf{f:alm_tr_prop}, therefore the set
$\CC(v_1, v_2)$ is dense in $C^0(M_1, M_2)$.

It is easy to show that if  $A$ satisfies  \rrf{f:tr_prop},
then for every $p\in S_k(f_1)$ and
$q\in S_{n_2-k}(f_2)$
the set
$T(p,q)=D(p, v_1)\cap A^{-1}(D(q,-v_2))$
is finite (where $n_2=\dim M_2$)
and therefore the intersection index
$$
n(p,q; A)=A(D(p,v_1))\krest D(q,-v_2)\in\ZZZ
$$
is defined.
Introduce  a  homomorphism of graded groups
\bq\lb{f:def_flat}
A_\flat:\MM_*(f_1, v_1)\to \MM_*(f_2, v_2)
\end{equation}
 by the formula
\bq
A_\flat(p)=\sum_{q}
n(p,q;A)\cdot q,
\end{equation}
where  the
summation in the formula is over
$q\in S_{n_2-\ind p}(f_2)$.

\bepr\mlb{p:sharp_flat}
$A_\sharp=A_\flat$.
\enpr
\Prf
Let $B=\Phi(T,-v_2)\circ A$; then
$B_\flat=A_\flat,
B_\sharp=A_\sharp$.
For sufficiently large $T$ the map $B$ satisfies
$B(M_1^{(k)})\sbs M_2^{(k)}$ for every $k$.
For $p\in S_k(f_1)$
let $\D_p\in H_k(M_1^{(k+1)}, M_1^{(k)})$
be the homology class of the descending disc of $p$.
It is easy to check the following formula
for the image of $\D_p$
\wrt~ $B_*$
\bq\lb{f:homol_inters}
B_*(\D_p)
=
\sum_{q} n(p,q;A) \D_q ~\in~ H_k(M_2^{(k+1)}, M_2^{(k)});
\end{equation}
here
$\D_q\in H_k(M_2^{(k+1)}, M_2^{(k)})$
is the homology class of the descending disc of $q$.
Thus by definition we have $B_\flat=B_\sharp$. $\qs$

%ii\input ncom

\section{Novikov complex}
\mlb{s:nov_com}

Let $M$ be a closed manifold, $f:M\to S^1$
a Morse map, and $v$ be an  oriented $f$-gradient.
Recall that if  $v$ satisfies \ta,
then the {\it Novikov complex} is defined
(see \cite{novidok}, \cite{noviuspe}, \cite{patou}).
This  is a chain complex of free finitely generated
modules
over
the ring
$$\ove L=\ZZZ((t))=\ZZZ[[t]][t^{-1}].$$
The definition of this complex
is based on a  procedure of  counting
the  flow lines of $v$ joining  critical points
of $f$. In   Subsection \mrf{su:constru_novik} we
 construct the Novikov complex
in a more general situation, when
$v$ is only almost transverse.
In  Subsection \mrf{su:canon_chain_equiva_nov}
we construct a canonical chain equivalence between the Novikov
complex and the completed singular chain complex of the
infinite cyclic covering corresponding to $f$.
In the third subsection we study
the functorial properties of the Novikov
complex thus obtained.

\subsection{Construction of the  complex }
\mlb{su:constru_novik}

Let $\bar M\to M$ be
the infinite cyclic covering, induced by
$f$ from the universal covering
$\RRR\to S^1$.
Lift the function $f$ to a Morse function
$F:\bar M\to \RRR$.
Let $t$ be the generator of the structure group
$\approx \ZZZ$ of the covering, such that
$F(tx)=F(x)-1$ for every $x\in \bar M$.
Let
$$
P=\ZZZ[t],~ \wh P=\ZZZ[[t]],~
L=\ZZZ[t, t^{-1}],~
P_n=P/t^nP
\mxx{ where } n\geq 0.
$$
Lift the vector field $v$ to a  vector field $\bar v$
on $\bar M$; then $\bar v$ is an \orial~
$F$-gradient, invariant \wrt~ the action of
the structure group of the covering.
Let
$\MM_k(F)$
be the free abelian group generated by the
critical points of $F$ of index $k$,
then $\MM_k(F)$
has a natural structure of a free $L$-module.
Let
$$
\NN_k(F)
=
\MM_k(F)\tens{L} \ove L.
$$
We shall define a boundary operator in
the graded group
$
\NN_*(F)$
thus endowing it with
the structure of a chain complex of $\ove L$-modules.
As a first step we shall  define
an {\it incidence coefficient}
$n(p,q;\bar v)$
for every pair $p,q$
of critical points of $F$ with $\ind p=\ind q+1$.
Choose any regular values $\l,\m$
of $F$
such that
$$
\l< F(p), F(q)<\m.
$$
Let $W=F^{-1}([\l,\m])$; then $\bar v|W$
is an \orial~ $F|W$-gradient, and the incidence coefficient
$n(p,q;\bar v|W)$ is defined
(see Definition \mrf{d:inc_coeff}).
It is clear that this integer
does not depend on the particular choice of
the regular values $\l,\m$;
it will be denoted $n(p,q;v)$.

\bepr\mlb{p:def_novik}
The formula
$$
\pr_k p=
\sum_q n(p,q;v) q
$$
defines a homomorphism
$\pr_k:\NN_k\to \NN_{k-1}$
of $\ove L$-modules
such that
$\pr_k\circ \pr_{k-1}=0$.
\enpr
Let $\l$ be a regular value of $F$,
denote by
$\NN_*^\l$
the set of all $\xi\in \NN_*$
such that $\supp\xi\sbs F^{-1}(]-\infty, \l])$.
It is clear that $\NN_*^\l$ is a
free and finitely generated
$\wh P$-submodule of $\NN_*$, which is $\pr_*$-invariant.
The module $\NN_*^\l/t^n \NN_*^\l$
is a free $P_n$-module, and by definition
the truncated boundary operator
$$
\pr_k^{(n)}
:
\NN_k^\l/t^n \NN_k^\l
\to
\NN_{k-1}^\l/t^n \NN_{k-1}^\l
$$
is equal to  the boundary operator in the Morse
complex of the restriction of $F$ and $\bar v$
to the cobordism $W_n=F^{-1}([\l-n,\l])$.
Therefore $\pr_{k-1}^{(n)}\circ \pr_{k}^{(n)}=0$,
and this implies the assertion. $\qs$
\bede\mlb{d:def_nov}
The graded $\ove L$-module $\NN_*$ endowed with
the boundary operator $\pr_*$
will be denoted $\NN_*(f,v)$
and called {\it Novikov complex}.
\end{defi}

\subsection{Canonical chain equivalence}
\mlb{su:canon_chain_equiva_nov}

In this subsection  we  construct
a chain equivalence
\bq\lb{f:cha_eq}
\NN_*(f,v)\to \SS_*(\bar M)\tens{L} \ove L,
\end{equation}
where $\SS_*(\bar M)$ is the singular chain complex
of $\bar M$.
This construction is  quite similar to the case of
transverse gradients, and we shall
give only an outline of the construction,
following the exposition in \cite{patou}, and \cite{paclo}.
The chain equivalence \rrf{f:cha_eq}
is  built from the
Morse equivalences
corresponding to the  cobordisms
$W_n=F^{-1}([\l-n, \l])$
(where $\l$ is any regular value of $F$).
The truncated Novikov complex
$\NN_*^\l/t^n \NN_k^\l$
is basis-preserving isomorphic to the
Morse complex $\MM_*(F|W_n, \bar v|W_n)$
(as chain complexes of $\ZZZ$-modules).
The composition of this isomorphism
with the Morse chain equivalence
will be denoted
$$
\ve_n:\NN_*^\l/t^n \NN_*^\l
\to
\SS_*(W_n, \pr_0 W_n).
$$
Let
$$
V_\l^-=F^{-1}(]-\infty, \l]).
$$
The composition
\bq\lb{f:composition}
\NN_*^\l/t^n \NN_*^\l
\rTo^{\ve_n}
\SS_*(W_n, \pr_0 W_n)
\rInto
\SS_*(V_\l^-, t^nV_\l^-)
\end{equation}
is a chain equivalence.
Both the source and the target of this chain equivalence
are free $P_n$-modules, but the  map
\rrf{f:composition}
 is apriori only a homomorphism
of abelian groups.
We shall now prove that the map
\rrf{f:composition}
is chain homotopic
to a chain equivalence
of $P_n$-modules.
To this end recall that the Morse chain equivalence
$\ve_n$ is defined using a Morse-Smale filtration
of $W_n$, adjusted to $(f,v)$.
Let $W_n^{(k)}$ be such a filtration. Define a filtration
$X^{(k)}$
of the space  $V_\l^-$
as follows:
$$
X^{(k)}
=
W_n^{(k)}
\cup
t^n V_\l^-.
$$
\bede\mlb{d:equiva_ms_filtr}
A Morse-Smale filtration $W_n^{(k)}$ of $W_n$
is called {\it equivariant} if
$$
tX^{(k)}\sbs X^{(k)}
\mxx{ for every } k.
$$
\end{defi}
The following proposition is proved in
\cite{paclo}  (Proposition 2.3) for the case when $v$
satisfies \ta; the proof is carried over
to the case of almost
transverse gradients without much changes.
\bepr\mlb{p:exist_equiv_filtr}
For every Morse-Smale filtration $W^{(k)}_n$
of the cobordism $W_n$
there is an equivariant
Morse-Smale filtration $\wi W^{(k)}_n$
such that $\wi W^{(k)}_n\sbs W^{(k)}_n$. $\qs$
\end{prop}

If   $W_n^{(k)}$
is an equivariant Morse-Smale
filtration
then the filtration
$$
\SS_*(X^{(k)},t^n V_\l^-)
$$
of the singular chain complex
$\SS_*(V_\l^-,t^n V_\l^-)$
is a filtration  by $P_n$-modules.
It is not difficult to show that this
filtration is a {\it cell-like } filtration
by $P_n$-modules
(that is, the homology
of the pair $(X^{(k+1)}, X^{(k)})$
vanishes in all degrees except $k$,
and in this degree
the homology is a free $P_n$-module).
Therefore the  adjoint chain complex
$H_*(X^{(*+1)},X^{(*)})$
is basis-preserving isomorphic to
$\NN_*^\l/t^n \NN_*^\l$, and
the canonical chain equivalence
\bq\lb{f:e_n}
\EE_n
:
\NN_*^\l/t^n \NN_*^\l
\to
\SS_*(V_\l^-, t^n V_\l^-)
\end{equation}
of $P_n$-modules associated with
this filtration
is homotopic over $\ZZZ$
to the chain equivalence
\rrf{f:composition}.
$\qs$

The proof of the next proposition
repeats the proof of
Theorem \mrf{t:independ_ord}
with corresponding modifications.
\bepr\mlb{p:e_n_independ}
The chain equivalence
\rrf{f:e_n}
is independent
(up to chain homotopy) of the particular choice
of the equivariant Morse-Smale filtration
adjusted to $(f|W_n,\bar v|W_n)$.
\end{prop}

\bede\mlb{d:equiv_equiv}
The chain equivalence \rrf{f:e_n}
will be called {\it equivariant Morse equivalence}.
\end{defi}
Now we shall compare the
equivariant Morse equivalences corresponding
to different regular values of $F$.
\bepr\mlb{p:compare_l_m}
Let $\l\leq \m$ be regular values of $F$
and $m\leq n$ be positive integers.
Then the following diagram of chain maps of
complexes over $P$ is
homotopy commutative:
\bq\lb{f:commut_equiv}
\xymatrix{
\NN_*^\l/t^n \NN_*^\l
\ar[r] \ar[d] &
\NN_*^\m/t^m \NN_*^\m \ar[d]\\
\SS_*(V_\l^-, t^n V_\l^-)
\ar[r] &
\SS_*(V_\m^-, t^m V_\m^-)}
\end{equation}
(here the vertical arrows stand for
the equivariant Morse equivalences,
and the horizontal arrows are induced by the inclusion
$V_\l^- \rInto V_\m^-$).
\enpr
\Prf See \cite{paclo}, Lemma 3.8. $\qs$

Now we can proceed to the description of the
canonical chain equivalence
between the Novikov complex and the
completed singular chain complex.
For $\l\in \RRR$
put
$$
\SS_*^\l=\SS_*(V_\l^-).
$$
Then  $\SS_*^\l$
is a chain complex of free $P$-modules.
Let
$$
\wh\SS_*^\l
=
\liminv~\SS_*^\l/t^n \SS_*^\l
$$
be the completion of this module, then
the natural inclusion
$$
\SS_*^\l\tens{P}\wh P\sbs \wh\SS_*^\l
$$
is a homotopy equivalence.
Note also that the module $\SS_*^\l\tens{P}\wh P$
is identified in a natural way with the
submodule of all the elements of
$\SS_*(\bar M)\tens{L}\ove L$
with support below $\l$.
\bede\mlb{d:canonical_equiv}
A chain map
$$
A_*:\NN_*(f,v)\to \SS_*(\bar M)\tens{L} \ove L
$$
is called {\it compatible with the
equivariant Morse equivalences}
if there is $\l$ such that
$$
A_*(\NN_*^\l)\sbs \SS_*^\l\tens{P}\wh P
$$
and for every positive integer $n$
the quotient chain map
$$
A_*/t^n :\NN_*^\l/t^n \NN_*^\l\to \SS_*^\l/t^n \SS_*^\l
$$
is homotopic over $P_n$
to the equivariant Morse equivalence $\EE_n$.
\end{defi}
\beth\mlb{t:canonic_equiv}
There  is a homotopy unique
chain equivalence
$$\gE_*:\NN_*(f,v)\to \SS_*(\bar M)\tens{L}\ove L$$
compatible with the equivariant Morse equivalences.
\end{theo}
\Prf
1. {\it Existence.}
Let $\l$ be any regular value of $F$.
As it follows from Proposition \mrf{p:compare_l_m}
for every positive $n$
the following diagram is
chain homotopy commutative
$$
\xymatrix{
\NN_*^\l/t^n \NN_*^\l
 \ar[d]^{\EE_n} &
\NN_*^\l/t^{n+1} \NN_*^\l \ar[l]\ar[d]^{\EE_{n+1}}\\
\SS_*(V_\l^-, t^n V_\l^-)
 &
\SS_*(V_\l^-, t^{n+1} V_\l^-) \ar[l]}
$$
Then one can prove  that there
is a chain map
$$
\wh\EE_*:\NN_*(f,v)=\liminv~\NN_*^\l/t^n \NN_*^\l
\to
\liminv~\SS_*(V_\l^-)/ t^n \SS_*(V_\l^-)
=
\wh \SS_*(V_\l^-)
$$
such that for every $n$
the quotient map $\wh\EE_*/t^n$ is
homotopic to $\EE_n$
(see  \cite{patou}, \S 3 B).
The composition
$$
\NN_*(f,v)=\NN_*^\l\tens{\wh P}\ove L
\rTo^{\wh \EE_*\otimes \id}
\wh\SS_*(V_\l^-)\tens{\wh P}\ove L
\rTo^\sim \SS_*(V_\l^-)\tens{P}\ove L²
\approx
\SS_*(\bar M)\tens{L}\ove L
$$
is then a chain equivalence
compatible with the equivariant Morse equivalences.

2. {\it Uniqueness}.
Let $\gE'_*:\NN_*\to \ove\SS_*$
be another chain map
compatible with the equivariant Morse
equivalences, and let
$\m$ be the corresponding
regular value of $F$. Let
$$
\AA=\gE_*~|~\NN_*^\l,
\quad
\AA'=\gE'_*~|~\NN_*^\m
$$
Assuming that $\l\leq \m$,
consider the following diagram:
\bq\lb{f:uniq_dia}\xymatrix{
\NN^\l_*/t^n\NN^\l_*
\ar[d]^{\AA/t^n}\ar[r]
&
\NN^\m_*/t^n\NN^\m_* \ar[d]^{\AA'/t^n}\\
 \SS_*^\l/t^n\SS_*^\l \ar[r]
&
\SS_*^\m/t^n\SS_*^\m }
\end{equation}
where the both horizontal
arrows are induced by the inclusion
$V_\l^-\sbs V^-_\m$.
The vertical arrows are by the hypotheses homotopic to
the equivariant
Morse equivalences, so  the  diagram is
identical with the diagram \rrf{f:commut_equiv}
(with $n=m$)
and therefore
\rrf{f:uniq_dia}
is homotopy commutative.
Therefore the compositions
$$\NN^\l_*\rTo^{\AA}\SS_*^\l\rInto
\ove\SS_*^\m$$
and
$$\NN^\l_*\rInto\NN^\m_*
\rTo^{\AA'}\ove\SS_*^\m$$
have the property that their $t^n$-quotients
are chain homotopic for every $n$.
The proof is then completed by the following lemma
(\cite{paclo}, Prop. 2.8).
\bele\mlb{l:noet}
Let $A_*, B_*$
be chain complexes of  $\wh P$-modules,
which are homotopy equivalent to
finitely generated free complexes.
Let $\psi, \psi':A_*\to B_*$
be chain maps. Assume that for every $n$
the quotient chain maps
$\psi/t^n, \psi'/t^n:A_*/t^n A_*\to B_*/t^n B_*$
are chain homotopic. Then $\psi$ and $\psi'$
are chain homotopic.\enle
$\qs$

\bere\mlb{r:schuetz}
Replacing  the singular chain complex
by the simplicial chain complex $\D_*(\bar M)$, associated to
a $\smo$ triangulation of $M$,
we obtain a chain equivalence
$$
\NN_*(f,v)\rTo^\sim \D_*(\bar M)\tens{L} \ove L.$$
In \cite{schuetzCI}
D.Sch\"utz gives an explicit formula
for the inverse of this  chain equivalence
(in the case when $v$ is transverse):
\bq\lb{f:schuetz_triang}
\xi(\s)=\sum_q (\s\krest D(q,-v))
\end{equation}
where $\s$ stands for a simplex of the
triangulation of $\bar M$,
$\dim\s=k$, the summation in
the right hand side is over all the
critical points of $F:\bar M\to \RRR$
of index $k$, and $\krest$
is the algebraic intersection number.
The formula \rrf{f:schuetz_triang} is valid
in the assumption
that the simplices of the triangulation
are  transverse to
the ascending discs of all critical points of $F$.
It is not difficult to show that the
formula \rrf{f:schuetz_triang}
remains valid if
the $f$-gradient $v$ is almost transverse.
\enre
\subsection{Functorial properties}
\mlb{su:functo_nov}

In this subsection we shall study
the functorial properties
of the Novikov complex and the
canonical chain equivalence $\gE_*$.
Let $f_1:M_1\to S^1, f_2:M_2\to S^1$
be circle-valued Morse maps (where $M_i$
are  closed manifolds).
Let $v_i$  be
\orial~ $f_i$-gradients (where $i=1, 2$).
Let $\smo_0(M_1, M_2)$
denote  the set of all $\smo$ maps $A$
such that
\bq\lb{f:cond_induce}
f_2\circ A \sim f_1.
\end{equation}
This is an open set in $\smo(M_1, M_2)$,
and we shall assume that it is not empty.
Let $A\in \smo_0(M_1, M_2)$.
We impose the following condition on $A$:
\begin{multline}\lb{f:tr_propp}
A|D(p,v_1):D(p,v_1)\to M_2
\mx{
is transverse to } D(q,-v_2)\\
\mx{
for every } p\in S(f_1), q\in S(f_2).
\end{multline}
A routine argument from
transversality theory proves the next propostion.
\bepr\mlb{p:dens_nov}
The set of the maps
$A\in \smo_0(M_1, M_2)$
satisfying \rrf{f:tr_propp}
is residual
in $\smo_0(M_1, M_2)$.$\qs$
\enpr
Choose and fix  a lift
$\bar A:\bar M_1\to \bar M_2$
of the map $A$ to the infinite cyclic coverings
(such lift exists since $f_2\circ A\sim f_1$).
Let $F_1, F_2$ be any lifts of
the functions $f_1$, resp. $f_2$ to
the coverings $\bar M_1, \bar M_2$.
Set $\dim M_1=n_1,~ \dim M_2=n_2$.
Let $p\in S_k(F_1), q\in S_{n_2-k}(F_2)$.
It is easy to see that
if the map $A$ satisfies \rrf{f:tr_propp}, then
the set
$$
T(p,q)=D(p,\bar v_1) \cap \bar A^{-1}\big( D(q,-\bar v_2) \big)
$$
is finite.
Therefore the algebraic intersection index
$$
N(p,q;\bar A)=\bar A\big(D(p,\bar v_1)\big)
\krest D(q, -\bar v_2)\in \ZZZ
$$
is defined. Observe  that
$N(t^k p, t^k q; \bar A) =N(p,q;\bar A)$.
For $p\in S_k(F_1)$ put
\bq\lb{f:defi_induced}
\bar A_\sharp(p)
=
\sum_{q} N(p,q;\bar A) q
\end{equation}
(where the summation is
over the set of all critical points of $F_2$
of index $k$).
It is clear that the expression
in the right hand side of \rrf{f:defi_induced}
is an element of
$\NN_k(f_2,v_2)$,
and the formula
\rrf{f:defi_induced}
defines  a homomorphism of graded
$\bar L$-modules
$$\bar A_\sharp:~\NN_*(f_1, v_1)\to \NN_*(f_2,v_2).$$
\beth\mlb{t:functo_nov}
The graded homomorphism $\bar A_\sharp$
is a chain map, and the following diagram is
chain homotopy commutative:
\bq\lb{f:functo_nov}
\xymatrix{
\NN_*(f_1, v_1) \ar[r]^{\bar A_\sharp}\ar[d]^{\gE_*^{(1)}}
& \NN_*(f_2, v_2) \ar[d]^{\gE_*^{(2)}}\\
\SS_*(\bar M_1)\tens{L}\ove L
 \ar[r]^{\bar A_*} & \SS_*(\bar M_2)\tens{L}\ove L
}
\end{equation}
(where $\gE_*^{(1)}$ and $\gE_*^{(2)}$
are the canonical chain equivalences).
\end{theo}
\Prf
Let $\l$ be a regular value  of $F_1$, and
$\mu$ be a regular value  of $F_2$.
We have the Novikov complexes
$$
\NN_*^{(1)}
=
\NN_*^\l(f_1, v_1), \quad
\NN_*^{(2)}
=
\NN_*^\m(f_2, v_2).
$$
Put
$$V^-_\l=F_1^{-1}(]-\infty, \l]),
\quad
U^-_\mu= F_2^{-1}(]-\infty, \mu]).$$
Choose   $\l$ and $\mu$ in such a way, that
$$
\bar A(V^-_\l)\sbs U^-_\mu,
\mxx{ so that }
\bar A_\sharp (\NN_*^{(1)})\sbs \NN_*^{(2)}.
$$
Consider the following diagram
(where $N$ is any positive integer):
\bq\lb{f:functo_nov_trunc}
\xymatrix{
\NN_*^{(1)}/t^N \NN_*^{(1)}
\ar[r]^{{\bar A_\sharp}/t^N}
\ar[d]^{\EE_N^{(1)}}
&
\NN_*^{(2)}/t^N \NN_*^{(2)}
\ar[d]^{\EE_N^{(2)}}\\
\SS_*(V^-_\l,V^-_{\l-N}) \ar[r]^{\bar A_*}
&
\SS_*(U^-_\m,U^-_{\m-N})
}
\end{equation}
(here $\EE_N^{(1)},~ \EE_N^{(2)}$
are the equivariant Morse equivalences).
To prove our theorem it suffices to show that
the upper horizontal arrow in this diagram is a
chain map, and that the diagram is chain homotopy commutative.
Pick any $t$-ordered Morse functions on the cobordisms
$$
W_{(1,N)}=F_1^{-1}([\l-N, \l]),\quad
 W_{(2,N)}=F_1^{-1}([\m-N, \m])$$
 and let $(X^{(k)}, V^-_{\l-N}),~ (Y^{(k)},U^-_{\m-N})$
 be the corresponding filtrations
 of the pairs $(V_\l^-,V^-_{\l-N})$, resp. $(U^-_\m,U^-_{\m-N})$.
Let us first consider a particular case
when the map $\bar A$ conserves these filtrations,
that is
\bq\lb{f:cons_filll}
\bar A(X^{(k)},V^-_{\l-N})
\sbs(Y^{(k)},U^-_{\m-N}).
\end{equation}
The adjoint complexes of these
filtrations are isomorphic
to the chain complex
$\NN_*^{(1)}/t^N \NN_*^{(1)}$
resp.
$\NN_*^{(2)}/t^N \NN_*^{(2)}$
via an isomorphism which preserves the bases.
The map $\bar A$ induces a chain map
$$A_\star:\NN_*^{(1)}/t^N \NN_*^{(1)}
\to
\NN_*^{(2)}/t^N \NN_*^{(2)}$$
of the adjoint complexes.
If we replace
the map $\bar A_\sharp/t^N$  by  $A_\star$  in
the diagram
\rrf{f:functo_nov_trunc}
the resulting diagram will be homotopy
commutative (by functoriality).
\bele\mlb{l:phi_phi}
$\bar A_\sharp/t^N=A_\star$.
\enle
\Prf
It suffices to check that
for every critical point $p$ of $F_1$
belonging to $V^-_\l\sm V^-_{\l-N}$
the fundamental class $\D(p,v_1)$
is sent to the same element by
$\bar A_\sharp$ and $A_\star$.
The map $A_\star$
sends this fundamental class
 to the fundamental class
of the manifold $\bar A(D(p,\bar v_1))$ in
the pair $(Y^{(k)},Y^{(k-1)})$.
The set $Z=Y^{(k)}\sm \Int Y^{(k-1)}$
is a cobordism, endowed with a Morse
function $F_2|Z$, which has only critical points
of indices $k$.
A standard computation shows  that
the $\bar A$-image of the fundamental class
of $D(p,\bar v_1)$
equals to
$$
\sum_q \Big( \bar A(D(p,\bar v_1))\krest \D(q,-v_2)\Big)
\D(q,v_2)
$$
where $q\in S_k(F_2)\cap Z$, and our
assertion  follows. $\qs$

Thus we have proved the homotopy commutativity
of the diagram
\rrf{f:functo_nov_trunc}
in the case when
\rrf{f:cons_filll}
holds.
Let us now proceed to the general case.
Let $B:M_1\to M_2$
be the map defined by
$$
B=\Phi(T,-v_2)\circ A,
$$
(where $T>0$).
Lift it to the map of coverings as follows:
$$
\bar B= \Phi(T,-\bar v_2)\circ\bar A
$$
Observe that the map
$ B$ satisfies the condition
\rrf{f:tr_propp}
and that
$$\bar A\big(D(p,-v_1)\big)
\krest D(q, -v_2)
=
\bar B\big(D(p,-v_1)\big)
\krest D(q, -v_2).
$$
If $T>0$ is sufficiently large,
the map $\bar B$
satisfies condition
\rrf{f:cons_filll}
and the proof of our theorem is now over. $\qs$

%ii\input further

\section{Applications and further developments}
\mlb{s:further}

Here we outline some applications
of the techniques developed in the preceding
sections.

\subsection{On the  $C^0$-stability in Novikov complex}
\mlb{su:cell_stab}

It is most likely that
the Novikov complex $\NN_*(f,v)$
is not stable \wrt~ the $C^0$-small perturbations of
$v$
(although I   know no explicit examples).
Still for every Morse function $f:M\to S^1$
there is a class of almost transverse gradients
for which the Novikov complex is $C^0$-stable;
namely the class of
gradients satisfying the condition $(\gC \CC)$,
introduced in \cite{pastpet}, \cite{paadv}.

Without reproducing here the definition
we shall just give a geometric  description of
this class. Let $f:M\to S^1$
be a Morse map, $v$ be an $f$-gradient.
Let $V=f^{-1}(\l)$
be a regular level surface of $f$.
The map which associates to $x\in V$
the point of the second  intersection with $V$ of the
$(-v)$-trajectory starting at $x$ is a
partially defined smooth map
$\stv : V\to V$. The  condition $(\gC \CC)$ requires
an existence of a certain
handle-like filtration $V^{(k)}$
of the manifold $V$
such that the map $\stv$
gives rise to a family of continuous maps
$V^{(k)}/V^{(k-1)}\to V^{(k)}/V^{(k-1)}$
between the successive quotients of the filtration.
See the  definition
in \cite{paadv}, p. 107, (or in \cite{pawitt}, page 317,
where the condition $(\gC \CC)$ was denoted $(\gC')$).
It is proved in \cite{paadv}, that the set of
$f$-gradients satisfying $(\gC \CC)$
is open and dense \wrt~ $C^0$-topology
in the set of all $\smo$ gradients.
It is also proved there (Theorem 5.7 of \cite{paadv})
that for every transverse
$f$-gradient
$v$ satisfying $(\gC \CC)$ and for every
open \nei~ $U$ of $S(f)$ there is $\d>0$
such that for every transverse
$f$-gradient $w$ with $||w-v||<\d$
the Novikov complexes $\NN_*(f,v)$ and $\NN_*(f,w)$
are basis-preserving isomorphic
(we imply here that $v$ and $w$ are similarly oriented).

It is easy to check that any gradient
$v$ satisfying $(\gC \CC)$
is almost transverse, so we can apply the
results of Subsection
\mrf{su:constru_novik}
to construct the Novikov complex
$\NN_*(f,v)$. Using the methods of the present work
(Subsection \mrf{su:gen_morcom}) it is
not difficult to extend the
Theorem 5.7 of \cite{paadv}
as follows:
%of $(\gC \CC)$-gradients:
\beth\mlb{t:nov_stable}
Let $f:M\to S^1$ be a Morse function,
let $v$ be an oriented  $f$-gradient
satisfying $(\gC \CC)$.
Then there is $\d>0$ such that for every
oriented
$f$-gradient $w$ with $||w-v||<\d$
the Novikov complexes
$\NN_*(f,v)$ and $\NN_*(f,w)$
are isomorphic via a basis-preserving isomorphism,
if $v$ and $w$ are similarly oriented.
\enth
The details of the proof of the theorem
 will be published elsewhere.

 \subsection{Lefschetz zeta functions
 of almost transverse gradients}
 \mlb{su:zeta}

 There is a remarkable
 relation between the simple homotopy type
 of  the Novikov complex and the
 Lefschetz zeta function, counting the closed orbits
 of the gradient flow. This
relation was discovered
by M.Hutchings and Y.J.Lee in the
work \cite{hulee}.
In this paper they  established
a formula which says  that
in the case when the
Novikov complex  is acyclic,
its Reidemeister
torsion equals the Lefschetz zeta function
of the gradient flow.
The author proved (\cite{pafest}) that
in the general case, when the
Novikov complex is not acyclic, the
Lefschetz zeta function
of the gradient flow
equals the Whitehead torsion of the canonical
chain equivalence between the Novikov complex and
the completed simplicial chain complex
of the infinite cyclic covering.
(In the paper \cite{pafest}
I considered  only the case of
$C^0$-generic gradients, and this restriction was
removed in the paper \cite{paclo},
using an approximation argument.)
These results were generalized to the case of
Lefschetz zeta functions
with values in the completions of group rings
by the author \cite{pawitt},
using the K-theoretic techniques developed by
A.Ranicki \cite{ranKNO}, and A.Ranicki and the author
\cite{pajandran}.
D.Sch\"utz (\cite{schuetzC}, \cite{schuetzCI})
generalized the results of \cite{pawitt}
 to the case of irrational forms
(see the papers \cite{schuetzC}, \cite{schuetzCI}),
using the techniques of 1-parameter fixed point theory
by R.Geoghegan and A.Nicas
(\cite{gnXCIV}, \cite{gnXCIVprim}).
In all these works only the case
of transverse gradients was considered.

In this subsection we outline a generalization
of these results  to the case of almost
transverse gradients.
Set $\Wh(\ove L)= \ove{K}_1(\ove L)/T$,
where $T$ is the
image in $\ove{K}_1(\ove L)$
of the group $\{\pm t^n|n\in \ZZZ\}\sbs \ove L^\bu$;
it is not difficult to check that
$\Wh(\ove L)$ is isomorphic to
the multiplicative group $\WW$
of all power series with integer coefficients
with the free term equal to $1$.
Let $f:M\to S^1$
be a circle-valued Morse function and $v$ be an \orial~
$f$-gradient.
Choose any $C^1$-triangulation of $M$ and lift it to
a $\ZZZ$-invariant triangulation of $\bar M$.
Let $\D_*(\bar M)$ be the corresponding
simplicial chain complex.
The composition of the   homotopy equivalence
$\gE_*$ from Theorem \mrf{t:canonic_equiv}
with the natural chain homotopy equivalence
$\SS_*(\bar M)\tens{L}\ove L \rTo^\sim
\D_*(\bar M)\tens{L}\ove L$
is a chain equivalence
$$
\gE'_*: \NN_*(f,v)\rTo^\sim  \D_*(\bar M)\tens{L}\ove L
$$
of two finitely generated based
chain complexes over
$\ove L$; the image of the Whitehead
 torsion of this chain equivalence
 in the group $\WW$ will be denoted
$w(f,v)$.

Now to zeta functions.
Choose a regular value $\l$
of the lift $F:\bar M\to\RRR$ of $f$.
The composition of gradient descent from $V_\l=F^{-1}(\l)$
to $V_{\l-1}$ with the map $t^{-1}:V_{\l-1}\to V_\l$
determines a (partially defined) diffeomorphism
$\Phi$ of $V_\l$ to itself.
Using the almost transversality
property it is not difficult to show that for any
$n\in\NNN$ the set of
fixed points of the $n$-th iteration
$\Phi^n:V_\l\to V_{\l-n}$
is compact. Let $L(\Phi^n)$  denote its index, and set
\bq\lb{f:def_z}
\z_L(-v)= \exp\Big(\sum_{n\geq 1}
\frac {L(\Phi^n)}n t^n\Big) \in \ZZZ[[t]].
\end{equation}
It is easy to check that the power series
$\zeta_L(-v)$ does not depend on the particular choice of
the regular value $\l$.

\beth\mlb{t:zeta_tors}
For every \orial~ $f$-gradient $v$ we have
$$
w(f,v)=(\zeta_L(-v))^{-1}.
$$
\enth

The proof of the theorem follows the lines of
the proof of the main theorem of \cite{paclo}: first
we check the theorem for  gradients
satisfying condition $(\gC\CC)$,
and then apply a perturbation argument.

%ii\input refcze

\end{document}

\end{document}